\documentclass[11pt]{article}
\usepackage{geometry}
 \usepackage{amsthm}
 \theoremstyle{plain}
\usepackage{ tikz, pgfplots, todonotes,color, amsmath}
\usetikzlibrary{decorations.markings,arrows,calc}
\definecolor{darkblue}{rgb}{0, 0, 0.3}

\usepackage{mathrsfs}
\usepackage{todonotes}
\usepackage{color}
\usepackage{caption}
\usepackage{amssymb, amsbsy}
\usepackage[colorlinks=true, pdfstartview=FitV, linkcolor=darkblue, citecolor=darkblue, urlcolor=darkblue]{hyperref}
\def\PEM{{\bf P}}

\def\CM{\mathbb G}
\def \tt{\widetilde \theta}
\def\ver{r}
\def\verB{R}
\def \o{\omega}
\def\Hcal{\mathcal H}
\def \bs{\boldsymbol}
\def\deltad{\delta}
\def \scr { \mathscr} 

\def\be{\begin{eqnarray}}
\def\ee{\end{eqnarray}}

\def\&{\hspace{-15pt}&}
\def\bea{\begin{eqnarray}}
\def\eea{\end{eqnarray}}
\def\Jum{J}
\def\s{\sigma}
\def\kat{\tilde{\kappa}}

\def \pa{\partial}
\def\C{{\mathbb C}}

\def\R{{\mathbb R}}
\def\N{{\mathbb N}}
\def\wh{\widehat}

\def\Z{{\mathbb Z}}

\def\a{\alpha}

\def\l{\lambda}

\def\1{{\bf 1}}
\def\Ec{{\mathcal E}}

\def\Mcal {{\mathcal M}}
\def \eqref #1{ (\ref{#1})}
\def \ds {\displaystyle}

\def\1{\mathbf 1}
\def\wt{\widetilde}
\def \wh {\widehat}

\def \Etavertex{ \eta_{_{\mathbf V }}}
\def \ddz{ \frac {\d z}{2i\pi}}
\def\g{\gamma}
\def \la{\label}

\def \QED {\hfill $\blacksquare$\par \vskip 3pt}
\def \G{\Phi}

\newtheorem{assumption}{Assumption}[section]

\def\tr {\mathrm {Tr}}
\def\M{\mathcal M}

\makeatletter
\@addtoreset{equation}{section}
\makeatother
\newtheorem{theorem}{Theorem}[section]
\newtheorem{theorem*}{Theorem}

\newtheorem{remark}{Remark}[section]

\newtheorem{proposition}{Proposition}[section]
\newtheorem{corollary}{Corollary}[section]
\newtheorem{definition}{Definition}[section]
\newtheorem{conjecture}{Conjecture}
\def\res{\mathop{\mathrm {res}}\limits_}
\def \m{\mathop}

\def\br{\begin{remark}\rm }
\def\er{\end{remark}}
\def\bd{\begin{definition}}
\def\ed{\end{definition}}
\def\bp{\begin{proposition}}
\def\ep{\end{proposition}}
\def\be{\begin{eqnarray}}
\def\ee{\end{eqnarray}}

\def\&{\hspace{-15pt}&}
\def\s{\sigma}

\def \pa{\partial}
\def\C{{\mathbb C}}

\def\R{{\mathbb R}}
\def\N{{\mathbb N}}
\def\wh{\widehat}
\def\Z{{\mathbb Z}}
\def\a{\alpha}
\def \d{\,\mathrm d}
\def\l{\lambda}
\def\1{{\bf 1}}
\def \pa{\partial}
\def \le{\left}

\def\ri{\right}

\def\Acal{{\mathcal A}}
\def\Fcal{{\mathcal F}}
\def\CP1{{\mathbb {P}^1}}
\def\Gcal{\mathcal G}
\def\ba{\begin{array}}
\def\ea{\end{array}}
\def\la{\label}
\def\p{\partial}
\def\phi{\varphi}

\makeatletter
\@addtoreset{equation}{section}
\makeatother

\def\f{\frac}

\def\a{\alpha}

\def\g{\gamma}
\def\s{\sigma}

\def\ka{\kappa}

\def\l{\lambda}
\def\At{\widetilde {\Acal}}
\def\Mt{\widetilde {\Mcal}}
\def\wt{\widetilde {\omega}}
\def \m{\mathop}
\def\a{\alpha}

\def\dim{{\rm dim}}
\def\tr{{\rm tr}}

\begin{document}
%\title{
%Tau-functions and 
%monodromy symplectomorphisms}

\vspace{0.2cm}
\begin{center}
\begin{huge}
{Tau-functions and 
monodromy symplectomorphisms}
\end{huge}\\
\bigskip
M. Bertola$^{\dagger\ddagger}$\footnote{Marco.Bertola@\{concordia.ca,sissa.it\}},  
D. Korotkin$^{\dagger}$ \footnote{Dmitry.Korotkin@concordia.ca},
\\
\bigskip
\begin{small}
$^{\dagger}$ {\it   Department of Mathematics and
Statistics, Concordia University\\ 1455 de Maisonneuve W., Montr\'eal, Qu\'ebec,
Canada H3G 1M8} \\
\smallskip
$^{\ddagger}$ {\it  SISSA/ISAS,  Area of Mathematics\\ via Bonomea 265, 34136 Trieste, Italy }\\
\end{small}
\vspace{0.5cm}
\end{center}

\author{M.Bertola and D.Korotkin}

\abstract{
We derive a new Hamiltonian formulation of  Schlesinger equations in terms of the dynamical $r$-matrix structure.
The corresponding symplectic form is shown to be  the pullback, under the monodromy map, of a natural symplectic form on the extended monodromy manifold.  We show that Fock-Goncharov coordinates are log-canonical for the symplectic form. Using these coordinates we define the symplectic potential on the monodromy manifold and interpret the Jimbo-Miwa-Ueno
tau-function as the generating function of the monodromy map. This, in particular, solves a recent conjecture by 
A. Its, O. Lisovyy and A. Prokhorov. }

\tableofcontents

\section{Introduction}

Symplectic aspects of the  monodromy map for the Fuchsian systems were studied starting from  \cite{ Hitchin,AlekMal2,KorSam};  in these papers it was proved that the monodromy map is a symplectomorphism from a  symplectic leaf in the space of coefficients of the system to a symplectic leaf in the monodromy manifold.
The    non-Fuchsian case was considered in 
\cite{FlaschkaNewell,Bondal,Ugaglia,Boalch2,Marta}. Remarkably, the simplest non-Fuchsian case of the Painlev\'e II hierarchy was treated in the paper \cite{FlaschkaNewell} in 1981, about 15 years before the Fuchsian case was studied in detail in \cite{ Hitchin,AlekMal2,KorSam}.
The key object associated to any Fuchsian or non-Fuchsian system of linear ODE's is the tau-function introduced by Jimbo, Miwa and Ueno \cite{JimboMiwa};
until now its significance in the framework of the monodromy symplectomorphism remained unclear; the main goal of this paper is to fill this gap
and prove a recent conjecture by Its-Lisovyy-Prokhorov \cite{ILP}. 

Our interest in this subject stems from  the study of  monodromy map of a second order equation  on a Riemann surface; such a map was also proven to be  a symplectomorphism 
\cite{Kawai,BKN,KorNeum,BK_TMF}. Remarkably,  the generating function of the  monodromy symplectomorphism (the "Yang-Yang" function) 
plays an important role in the theory of supersymmetric Yang-Mills equations  \cite{NRS} ; several steps towards understanding of this generating function were made in 
\cite{BKN}.

In this paper we address the question about the role of the generating function of the monodromy symplectomorphism in the context of Fuchsian equations 
on the Riemann sphere. The conclusion we arrive to is somewhat unexpected: such generating function can be naturally identified with the isomonodromic tau-function;
moreover, this interpretation allows to define the dependence of the tau-function on monodromy data.

The version of the monodromy map for Fuchsian systems used in the current paper is 
  slightly different from the monodromy map 
considered in \cite{Hitchin,AlekMal2,KorSam}. This version is
standard in the theory of isomonodromy deformations 
\cite{JimboMiwa} and it was also considered in  \cite{Jeffrey, Boalch2} from the symplectic point of view. 

\begin{figure}[htb]
\centering
 \includegraphics [scale=1.0]{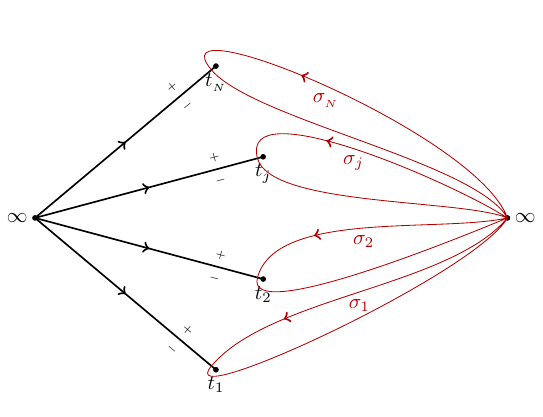}
 \caption{Choice of cuts and generators of the fundamental group.}
\label{genfun}
 \end{figure}

To describe the monodromy map we remind  the basics of the theory of solutions of Fuchsian systems of differential equations on
$\CP1$, following  \cite{JimboMiwa}.
Consider the  equation 
\be
\f{\p\Psi}{\p z}= \sum_{i=1}^N \f{A_i}{z-t_i}
 \Psi\;,\hskip0.7cm \Psi(z=\infty)=\1,
\la{lsint}
\ee
where $A_i\in sl(n)$ and $t_j\neq t_k$ such that $\sum_{i=1}^N A_i=0$.
Assume also that eigenvalues of each $A_j$ are simple and furthermore  do not differ by an integer.
Choose a system of cuts $\gamma_1,\dots,\gamma_N$  connecting  $\infty$ with $t_1,\dots,t_N$ respectively, and assume that the ends of these cuts
emanating from $\infty$ are ordered as $(1,\dots,N)$ counter-clockwise (Fig.\ref{genfun}).    
The normalization condition $\Psi(\infty)=\1$ is then understood in the sense that $\lim_{z\to\infty} \Psi(z)=\1$ where the limit is taken in the sector between  $\gamma_1$ and $\gamma_N$.

The set of  generators $\sigma_1,\dots,\sigma_N$ of the fundamental group
 $\pi_1(\CP1\setminus\{t_j\}_{j=1}^N,\infty)$ is chosen such that the loop representing $\sigma_j$ 
 crosses only the cut $\gamma_j$, and its orientation is chosen so that the relation 
 between $\sigma_j$ takes the  form $
 \sigma_N\cdot\dots\cdot \sigma_1={\rm Id}$ (Fig.\ref{genfun}).

The solution $\Psi$ of (\ref{lsint}) is single-valued in the simply connected domain
$\CP1\setminus\{\gamma_j\}_{j=1}^N$. Denote the diagonal form of the matrix  $A_j$ by $L_j$, $j=1,\dots, N$ (the matrices $A_j$ are diagonalizable due to our assumption about 
their eigenvalues). Then the asymptotics of $\Psi$ near $t_j$ has the standard form \cite{JimboMiwa}:
\be
\Psi(z)=(G_j+O(z-t_j)) (z-t_j)^{L_j} C^{-1}_j\;.
\la{asint}
\ee
The matrix $G_j$ is a diagonalizing matrix for $A_j$: 
\be
A_j=G_j L_j G_j^{-1}\;.
\la{AGint}
\ee
The matrices $C_j$ are called the {\it connection matrices}.
Notice that the matrices $G_j$ and $C_j$ are not uniquely defined by equation (\ref{lsint}) since a simultaneous transformation $G_j\to G_j D_j$
and $C_j\to C_j D_j$ with diagonal $D_j$'s  changes neither the asymptotics
(\ref{asint}) nor the equation (\ref{lsint}).

Analytic continuation of $\Psi(z)$ along $\sigma_j$  yields $\Psi(z) M_j^{-1}$, where the {\it monodromy matrix} $M_j\in SL(n)$ is related to the connection matrix $C_j$ and the exponent of monodromy $L_j$ by:
\be
M_j=C_j \Lambda_j C_j^{-1}\ ,\qquad \Lambda_j:= {\rm e}^{2i\pi L_j}.
\la{Mon}
\ee
Alternatively, the matrix $M_j$ can be viewed as the jump matrix on $\gamma_j$: orienting $\gamma_j$ from
$\infty$ towards $t_j$, the boundary values of $\Psi$ on the right ("+") and left ("-") sides of $\gamma_j$ are related by $\Psi_+=\Psi_i M_j$.
Our assumption about the ordering of the branch cuts $\gamma_j$ and generators $\sigma_j$ implies the relation
\be
M_1\cdots M_N=\1\;.
\la{prodmon}
\ee
The  {\it monodromy map}  introduced in \cite{JimboMiwa} sends the set of pairs
$(G_j,L_j)$ to the set of pairs $(C_j,\Lambda_j)$ for a given set   of poles $t_j$.

The map between the set of coefficients $A_j$ and the set of monodromy matrices $M_j$ is a different 
version of monodromy map associated to
equation (\ref{lsint});  the symplectic aspects of this version of the monodromy map were studied
in  \cite{Hitchin,AlekMal2,KorSam}.

To describe our framework in more details we introduce the following two spaces. The first space is the quotient 
\be
\Acal=\bigg\{(G_j,L_j)_{j=1}^N,\;\ \ G_j\in SL(n),  \ L_j \in \mathfrak h_{ss}^{nr}\ ,\forall j=1,\dots ,N:\ \ \ \sum_{j=1}^N G_j L_j G_j^{-1}=0\bigg\}\big/\sim
\la{defA}
\ee
where $\mathfrak h_{ss}^{nr}$ denotes the set of diagonal matrices with simple eigenvalues not differing by integers ({\it non-resonant}). The equivalence relation is given by the  $SL(n)$ action $G_j\mapsto S G_j$ with $S$ independent of $j$. 

The second space is the quotient
\be
\Mcal=\bigg\{\{C_j\,,\,L_j\}_{j=1}^N,\;\  C_j\in SL(n),\ L_j \in \mathfrak h_{ss}^{nr} :\ \prod_{j=1}^N C_j e^{2\pi i L_j} C_j^{-1}=\1\bigg\}\big/\sim\;.
\la{defM}
\ee
Similarly to (\ref{defA}), the equivalence is given by the $SL(n)$ action $C_j\mapsto S C_j$ (with the same $S$ for all $j$'s). 

For a fixed set of poles  $\{t_j\}_{j=1}^N$ we denote the "monodromy map" from the $(G,L)$-space to $(C,L)$-space by   
\be
\Fcal^t:\Acal\to \Mcal\;.
\la{Fdef}
\ee

\paragraph{\bf Poisson and symplectic structures on $\Acal$ and dynamical $r$-matrix.}

Let $\Hcal := SL(n,\C)\times \mathfrak h_{ss} = \{(G,L)\}$. Here $G\in {SL(N)}$ and $L={\rm diag}(\l_1,\dots,\l_n)$  is a diagonal traceless  matrix with $\l_j\neq \l_k$ and $\l_1+\dots+\l_n=0$.
Consider the following $1$-form on $\Hcal$:
\be
\theta=\tr (L G^{-1} \d G)\;.
\la{spGL}
\ee
We prove in Prop. \ref{nondeg} that the form $\o = \d \theta$ 
is non-degenerate, and therefore, is a symplectic form on $\Hcal$.
For any matrix $M$ we use the following notation for the Kronecker products
$$\m{M}^1 = M\otimes \1\;,\hskip0.7cm\m{M}^2 =\1\otimes  M.$$ 
Then the Poisson structure on  $\Hcal$ associated to the symplectic form $\o$ is  (see Prop. \ref{extPBKK}):
\be
\big\{\m{G}^1,\m{G}^2\big\}=-\m{G}^1\m{G}^2 r(L)\;,\hskip0.7cm
\big\{\m{G}^1,\m{L}^2\big\}=-\m{G}^1 \Omega\;,
\la{b2}
\ee
where 
\be
r(L)=\sum_{i<j} \f{E_{ij}\otimes E_{ji}-E_{ji}\otimes E_{ij}}{\l_i-\l_j} 
\la{dynr}\ee
and  
$$\Omega=\sum_{i=1}^n E_{ii}\otimes E_{ii} - \frac 1 n \1 \otimes \1\;;$$ 
we use the standard notation $E_{ij}$ for the  matrix with only one non-vanishing element  equal to $1$ 
in  the $(i,j)$ entry. The matrix $r(L)$ is a simplest example of dynamical $r$-matrix \cite{Etingof}.
Theorem \ref{thmrtoKK} shows that the bracket (\ref{b2})  induces the 
Kirillov-Kostant Poisson bracket for $A=GLG^{-1}$.

The bracket (\ref{b2}) can be used to define  the Poisson structure on the space $\Acal$ as follows. 
Denote first by $\Acal_0$  the space of pairs $\{(G_j,L_j)\}_{j=1}^N$ with the product symplectic structure, or, equivalently, with the following Poisson bracket:
\be
\{\m{G_j}^1,\m{G_k}^2\}=-\m{G_j}^1\m{G_k}^2 r(L_k)\,\delta_{jk}\;,\hskip0.7cm
\{\m{G_j}^1,\m{L_k}^2\}=-\m{G_k}^1 \Omega\,\delta_{jk}\;.
\la{b3}
\ee
The moment map corresponding to the group action $G_j\to SG_j$  ($S\in SL(n)$) on $\Acal_0$ is given by
$\sum_{j=1}^N G_j L_j G_j^{-1}$.
The space $\Acal$  (\ref{defA}) inherits a symplectic form from $\Acal_0$ via the standard symplectic reduction \cite{Audin}:
\begin{theorem*}
[see Thm. \ref{nondegmain}]
The Poisson structure induced on $\Acal$ from the Poisson structure (\ref{b3}) on $\Acal_0$  via the  reduction 
on the level set $\sum_{j=1}^N G_j L_j G_j^{-1}=0$ of the moment map, corresponding to the group action $G_j\to SG_j$,
is non-degenerate and the corresponding symplectic form is  
$
\o_{\Acal}=\d \theta_\Acal
$, 
where the  symplectic potential $ \theta_\Acal$ for  $\o_{\Acal}$  is given by 
$$
  \theta_\Acal=\sum_{k=1}^N\tr (L_k G_k^{-1} \d G_k)\;.
 $$
 \end{theorem*}
 This symplectic structure appeared in \cite{Boalch2} but the connection to dynamical $r$-matrix and associated Poisson structure  was not known until now.

\paragraph{\bf Symplectic structure on $\Mcal$.} Define the following 2-form  on the space $\Mcal$:
$$
\o_{\Mcal}=
\f{1}{4\pi i}(\o_1+\o_2)
$$
where
$$
\omega_1= \tr\sum_{\ell=1}^{N} 
\left( M_\ell^{-1} \d M _{\ell} \wedge K_{\ell }^{-1} \d K_{\ell} \right) 
+\tr  \sum_{\ell=1}^N  \left( \Lambda^{-1}_\ell C_{\ell} ^{-1}\d C_{\ell}\wedge  \Lambda_\ell  C_\ell ^{-1}\d C_\ell \right)
$$
$$
\o_2=2\sum_{\ell=1}^N  \tr  \left(\Lambda_\ell^{-1} \d \Lambda_\ell \wedge C_\ell^{-1} \d C_\ell
\right)
$$
with  $K_{\ell}=M_1\cdots M_\ell$ and $\Lambda_j=e^{2\pi i L_j}$.

The form $-\o_1/2$ coincides with the symplectic
form on the symplectic leaves $\Lambda_j=const$ of the  $SL(n)$ Goldman bracket 
(see (3.14) of \cite{AlekMal} in the case $g=0$).
The first  result of this paper  (see Thm.\ref{dTTH} and its proof in Section \ref{AMMalg}) is that
given a set of  poles $\{t_j\}_{j=1}^N$  and a point $p_0\in  \Mcal$ in a neighbourhood of which the monodromy map is invertible, 
the pullback of the form $\o_\Mcal$ under the map $\Fcal^t :\Acal\to \Mcal$ coincides with $\o_\Acal$:
\be
(\Fcal^t)^* \o_\Mcal=  \o_\Acal\;.
\la{symint}
\ee
This statement
 implies that (see Corollary \ref{cornond} and its proof)
the form $\o_\Mcal$ is closed and non-degenerate, and, therefore, defines a symplectic structure on $\Mcal$.

For a given set of monodromy data the monodromy map is invertible outside of a locus of codimension 1
in the space of poles \cite{Bolibruch}. Since  the form $\omega_{\Mcal}$ is independent of
$\{t_j\}$, this form is always  non-degenerate on the monodromy manifold.

The equality  (\ref{symint}) generalizes the results of 
\cite{Hitchin,AlekMal2,KorSam}, where it was proved that the monodromy map between the ``smaller" spaces - the space of coefficients $A_j$ with fixed eigenvalues and the symplectic leaf of the $GL(n)$  character variety of $N$-punctured sphere - is a symplectomorphism;
the formula (\ref{symint}) was proved in a different way in \cite{Boalch2}.

\paragraph{Time dependence.}
To describe  the  dependence  on the  $t_j$'s (the ``times'') we extend the spaces $\Acal$ and $\Mcal$ to include  also the coordinates $\{t_j\}$:
 \be
\widetilde {\Acal}=\Big\{(p,\{t_j\}_{j=1}^N)\;,\; p\in \Acal,\; t_j\in \C,\; t_j\neq t_k\Big\}\;,
\la{At}
\ee 
\be
\widetilde {\Mcal}=\Big \{(p,\{t_j\}_{j=1}^N)\;,\; p\in \Mcal,\; t_j\in \C,\; t_j\neq t_k\Big\}\;.
\la{Mt}
\ee

The monodromy map $\Fcal^t$ then  extends to the map 
$
\Fcal:\widetilde {\Acal}\to \widetilde {\Mcal}\;.
$
The locus in $\widetilde \Mcal$ where the map is not invertible is usually referred to as the {\it Malgrange divisor}.
Denote the pullback of the form $\o_\Acal$   from $\Acal$ to $\At$
by $\wt_\Acal$ and the pullback of the form $\o_\Mcal$  from $\Mcal$ to $\Mt$
by $\wt_\Mcal$ (notice that the forms $\wt_\Acal$ and $\wt_\Mcal$ are  closed but degenerate). Now we are in a position to formulate the next theorem (see Thm.\ref{ththa})
\begin{theorem*}[see Thm.\ref{ththa} together with Thm.\ref{dTTH}]
\la{conILP}
The following identity holds between two-forms on $\At$
\be
\Fcal^* \wt_\Mcal= \wt_\Acal-\sum_{k=1}^N d H_k\wedge dt_k
\la{wwint}
\ee
where
\be
H_k=\sum_{j\neq k} \frac{\tr A_j A_k}{t_k-t_j}\;,\hskip0.7cm k=1,\dots,N
\la{Hkint}
\ee
are the canonical Hamiltonians of the Schlesinger system.
\end{theorem*}
We remind that the Schlesinger equations \cite{Bolibruch}  consist of the following system of PDEs for the coefficients of $A(z)$
\be
\frac {\pa A_k}{\pa t_j} = \frac {[A_k,A_j]}{t_k-t_j},\ \ j\neq k\,;\hskip0.8cm  \frac {\pa A_j}{\pa t_j} = -\sum_{k\neq j}\frac {[A_k,A_j]}{t_k-t_j}
\la{Schlsys}
\ee
and they define the deformations of the connection $A(z)$ which preserve the monodromy representation. They are Hamiltonian equations with respect to the standard Kirillov-Kostant Poisson bracket with time--dependent Hamiltonians $H_k$ as in \eqref{Hkint}.

\paragraph{Tau function and generating function of the monodromy map.}
The above theorem allows to establish the relationship between the isomonodromic tau-function and the 
generating function of the monodromy map.
Namely, consider  some local  symplectic potential $\theta_\Mcal$ for the form $\o_\Mcal$ such that
$$
\d\theta_\Mcal=\o_\Mcal
$$  on the space 
$\Mcal$  (globally $\theta_\Mcal$ can be defined on a covering of
$\Mcal$) and denote its pullback to $\Mt$ by $\tt_\Mcal$.
Denote by $\tt_\Acal$ the pullback of $\theta_\Acal$ under the natural projection $\widetilde \Acal \to \Acal$. Then  $\tt_\Acal$ is the potential of  the symplectic form $\wt_\Acal$ on $\At$, 
and (\ref{wwint}) implies existence of a locally defined  generating function  $\Gcal$ on $\At$.

\begin{definition}
 The generating function (corresponding   to a  given choice of the symplectic potential $\theta_\Mcal$)  of the monodromy map  between spaces $\At$ and $\Mt$ 
is defined by
\be
d\Gcal=\tt_{\Acal}-\sum_{j=1}^N H_k dt_k -\Fcal^* \tt_\Mcal \;.
\la{defgen}
\ee
\end{definition}
A different choice of $\theta_\Mcal$ (and hence of its pullback $\widetilde \theta_\Mcal$) adds 
a $\{t_k\}$-independent term to $\Gcal$ i.e. it corresponds to a
transformation $\Gcal\to \Gcal+f(\{C,L\})$ for some local function $f$ on  $\Mcal$.

The dependence of $\Gcal$ on $\{t_j\}$ is, however, completely fixed by (\ref{defgen}). Namely, locally one can write (\ref{defgen}) in the coordinate system where $\{t_j\}_{j=1}^N$ and $\{C_j,L_j\}_{j=1}^N$ are considered as independent variables. Then derivatives of $G_j$ on $\{t_k\}$ for constant 
$\{C_j,L_j\}_{j=1}^N$ i.e. for constant monodromy data, are given by Schlesinger equations of isomonodromic deformations in $G$-variables:
\be
 \f{\p G_k}{\p t_j}= \f{A_j G_k}{t_j-t_k}\;,\;\;\; j\neq k\,;\hskip0.8cm
 \f{\p G_k}{\p t_k}=-\sum_{k\neq j} \f{A_j G_k}{t_j-t_k}\;.
 \la{eqGint}
 \ee
 The equations \eqref{eqGint} imply \eqref{Schlsys} {\it but not viceversa}. 
 In this paper we obtain the following Hamiltonian formulation of  equations (\ref{eqGint}):
 \begin{theorem*}[see Theorem \ref{Sheh}]
The equations \eqref{eqGint} are Hamiltonian,
 \be
 \frac{\pa G_k}{\pa t_j} = \{H_j, G_k \}  
 \ee
where $\{.,.\}$ is the quadratic   Poisson bracket \eqref{b2} and the Hamiltonians are given by \eqref{Hkint}. 
\end{theorem*}
A direct computation shows  that in $(t_j,C_j,L_j)$ coordinates the part of $\tt_\Acal$ containing $\d t_j$'s is given by $2\sum_{j=1}^N H_j \d t_j$; together with (\ref{defgen}) this implies
\be
\f{\p\Gcal}{\p t_j} = H_j\;.
\la{Gt}
\ee
Therefore,  we get the following theorem:
\begin{theorem*}[see Thm.\ref{tauGmain}]\la{tauG}
For any choice of symplectic potential $\theta_\Mcal$ on $\Mcal$ the dependence  of the  generating function   $\Gcal$ (\ref{defgen}) on $\{t_j\}_{j=1}^N$ coincides with the $t_j$-dependence of the 
isomonodromic Jimbo-Miwa tau-function. In other words, $e^{-\Gcal}\tau_{JM}$ 
depends only on monodromy data $\{C_j,L_j\}_{j=1}^N$. 
\end{theorem*}

The above theorem shows that the generating function  $\Gcal$  can be  used to define the  tau-function not only as a function of positions of singularities of the Fuchsian differential equation but also as a function of monodromy matrices.
The ambiguity built into this definition corresponds to the freedom to choose different symplectic potentials on  different open sets of  the monodromy manifold.

The symplectic potential we use in this paper was found in \cite{BK} using the coordinates
introduced by Fock and Goncharov in \cite{FG} (for $SL(2,\R)$ case these coordinates 
called {\it shear coordinates} are attributed to Thurston, see
\cite{Fock,CheFo}; see also \cite{Marta1} where the complex analogs of the shear coordinates were used for the explicit parametrization of  the open subset of full dimension of the  $SL(2,\C)$ character variety of four-punctured sphere).

\begin{definition}\la{deftauintro}
The $SL(n)$  tau-function $\tau$ on $\widetilde {\Mcal}$ is locally defined by the following set of compatible equations. The equations with respect to $t_j$ are given 
by the  formulas
\be
\f{\p \log \tau}{\p t_j}=\f{1}{2}\res{z=t_j}{\rm tr} A^2(z)\;,
\la{tautint1}
\ee
where
\be
A(z):= \sum_{i=1}^N \f{A_i}{z-t_i}.
\ee
The equations with respect to coordinates on monodromy manifold $\Mcal$ 
are given by
\be
\d_\Mcal \log \tau=\sum_{j=1}^N \tr (L_j G_j^{-1} \d_\Mcal G_j)-\theta_\Mcal[\Sigma_0]
\la{tauMcint1}
\ee
where $\theta_\Mcal[\Sigma_0]$ is a symplectic potential (\ref{symmpoten}) for the form $\omega_\Mcal$ 
defined using the Fock-Goncharov coordinates corresponding to a ciliated triangulation $\Sigma_0$ (see Section \ref{SecN}).
\end{definition}

Explicit formulas for derivatives of $\tau$ with respect to Fock-Goncharov coordinates will be given in Section \ref{tau_fun_mon}.

\paragraph{\bf Conjecture by A.Its, O.Lisovyy and A.Prokhorov.}
Theorem \ref{conILP} emphasizes  a close relationship  with  the  recent work \cite{ILP} where the issue of dependence of the Jimbo-Miwa tau-function on monodromy matrices was also addressed. In particular, the relevance of the Goldman bracket and the corresponding symplectic form on its symplectic leaves was observed in \cite{ILP} in the case of $2\times 2$ system with four simple poles (the associate isomonodromic deformations give Painlev\'e 6 equation).  

Moreover, the authors of \cite{ILP} introduced a form which we denote by $\Theta_{ILP}$  (this form is denoted by $\omega$ in  (2.7) of \cite{ILP}).
This form appeared in \cite{ILP} as a result of computation involving the 1-form introduced by  Malgrange
in \cite{Malg}, similarly to this work,  which in our notations is given by 
\be
\Theta_{ILP}=\sum_{j<k}^N\tr A_j A_k \d\log (t_j-t_k)+\sum_{j=1}^N \tr (L_j G_j^{-1} \d_{\Mcal} G_j)
\la{ILP}
\ee
where $\d_{\Mcal}$ denotes the differential with respect to monodromy data.
Proposition 2.3 of \cite{ILP} shows that the external derivative of the form (\ref{ILP})  is a closed 2-form independent of $\{t_j\}_{j=1}^N$.
Furthermore, in Section 1.6 the authors of \cite{ILP} formulate the following 
\begin{conjecture}
\la{ILPconj}
[Its-Lisovyy-Prokhorov] The form $\d\Theta_{ILP}$ coincides with the natural
symplectic form on the monodromy manifold.
\end{conjecture}

There are  two natural versions of this conjecture:
\begin{itemize}
\item
{\it The "weak" ILP conjecture.} In this version  $\d_\Mcal$ 
means the differential on a symplectic leaf $\{\Lambda_j=\mathrm{const}\}_{j=1}^N$ of the $SL(n)$ character variety of $\pi_1(\CP1\setminus\{t_j\}_{j=1}^N)$ (we denote this symplectic leaf by $\Mcal_\Lambda$). The canonical symplectic form on $\Mcal_\Lambda$ is given by the inversion of
the $SL(n)$ Goldman's bracket \cite{Goldman} and can be written explicitly 
in terms of monodromy data as shown in (\cite{AlekMal}, formula (3.14) for $g=0$ and $k=2\pi$).

By  {\it``weak" ILP conjecture} we understand the    coincidence of $\d\Theta_{ILP}$ (\ref{ILP})  with the
Goldman's symplectic form on  the symplectic leaves.

The problem with this formulation is that the choice of matrices $G_j$ should be such that they satisfy the 
Schlesinger equations (\ref{eqGint}); this requirement is not natural from the symplectic point of view.

\item
{\it The ``strong" ILP conjecture.} In this version the differential $\d_{\Mcal}$ in (\ref{ILP}) 
means the differential on the full space $\Mcal$ (\ref{defA}) which contains both the eigenvalues of the monodromy matrices  and the connection matrices. Then (omitting the pullbacks) the {\it strong ILP conjecture} states that
\be
\d\Theta_{ILP}=\o_\Mcal\;.
\la{strongILPf}
\ee
\end{itemize}

The weak version of the ILP conjecture can be  derived   from known results of \cite{Hitchin,AlekMal2} or \cite{KorSam}, as shown in Section \ref{previous}.

The strong version of the ILP conjecture is equivalent to our Theorem \ref{conILP}. 
To see this equivalence it is sufficient to write (\ref{defgen}) in coordinates which are split into "times" $\{t_j\}$ and some coordinates  on the monodromy manifold $\Mcal$. Then the "$t$-part" of the form $\widetilde {\theta}_\Acal$ is given by $2\sum_{k=1}^N H_k d t_k$ (this follows from the isomonodromic equations (\ref{eqGint}) for $\{G_j\}$) and the monodromy part coincides with the second term of the  form (\ref{ILP}) 
where the differential $\d_{\Mcal}$ is understood as the differential on $\Mcal$.
Now, taking the external derivative of (\ref{defgen}) we come to (\ref{wwint}) where the right-hand side coincides with the form $\d\Theta_{ILP}$ of \cite{ILP}. Finally, we notice that the formula (\ref{defgen}) allows to interpret the generating function 
$\Gcal$ as the action of the multi-time hamiltonian system, according to Conjecture 2 of \cite{IP} (see also \cite{IP1}).

Summarizing, the main results of this paper are the following:
\begin{enumerate}
\item
We give a new hamiltonian formulation of Schlesinger system written in terms of $(G,L)$-variables;
this formulation 
involves a quadratic Poisson structure defined by the dynamical $r$-matrix (Sec.\ref{DynSec}).
\item
We prove that the monodromy map for a Fuchsian system is a symplectomorphism between $(G,L)$ and $(C,\Lambda)$
spaces (Sec.\ref{AMMalg}). 
\item
We prove the ``weak" (Sec.\ref{previous}) and ``strong" (Thm.\ref{dTTH}) versions of the Its-Lisovyy-Prokhorov conjecture about coincidence of the
external derivative of the Malgrange form with the natural symplectic form on the monodromy manifold
\item
We introduce defining equations for the Jimbo-Miwa-Ueno tau-function with respect to Fock-Goncharov coordinates on the
monodromy manifold (Def.\ref{Deftau} and formula (\ref{tauNew})).
\item In the $SL(2)$ case we  derive equations which define the monodromy dependence of $\Psi$, $G_j$ and tau-functions (Thm.\ref{Psizeta}, Cor.\ref{Gzeta} and Prop.\ref{tau1zeta}).
\end{enumerate}

  \section{Dynamical $r$-matrix formulation of the Schlesinger system}
  \la{DynSec}
  In this section we describe the Hamiltonian formulation of   Schlesinger equations.
We start from considering the $GL(n)$ case and then indicate the modifications required in the  $SL(n)$ case.

  \subsection{Quadratic Poisson bracket via dynamical $r$-matrix }

Let us  introduce the space 
\be
\Hcal:= GL(n,\C) \times \mathfrak h_{ss} 
\ee
where $\mathfrak h_{ss}$ is the space  of diagonal matrices with distinct eigenvalues.
We denote an  element of $\Hcal$ by $(G,L)$ where $G\in GL(n)$ and $L\in h_{ss}$.
\bp
\la{nondeg}
Consider the following one-form on $\Hcal$:
\be
\theta := \tr (L G^{-1} \d G)\;.
\ee
Then the 2-form $\o=\d\theta$ given by
\be
\omega = \tr (  \d L \wedge G^{-1} \d G ) -\tr (L G^{-1} \d G \wedge G^{-1} \d G )
\la{defoGL}
\ee
is symplectic  on $\Hcal$.
\ep
{\bf Proof.}
The form $\o$ is obviously closed; to verify its non-degeneracy we consider two tangent vectors in $T_{(G,L)}\Hcal$ and represent them as 
$(X_i,D_i)\in gl(n) \oplus \mathfrak h$ ($i=1,2$) where $ \mathfrak h$ denotes the Cartan subalgebra of $gl(n)$. 
Then
\be
\o( (X_1,D_1), (X_2, D_2) )= \tr \bigg(D_1 X_2 -D_2 X_1 - L[X_1,X_2]\bigg)=\tr \bigg(D_1 X_2 - \big( D_2  + [X_2,L]\big)X_1\bigg)\;.
\la{pariin}
\ee
Suppose  that $\o$ is degenerate i.e. the vector $(X_2, D_2)$  can be  chosen so that (\ref{pariin})  vanishes identically for all $(X_1,D_1)$. Then, choosing $D_1=0$, we have  $\tr (( D_2 + [X_2,L]) X_1)=0$. Then, since $X_1$ is arbitrary, we have $ D_2 + [X_2,L]=0$; since $L$ is diagonal, the commutator is diagonal-free and hence $ D_2 =0$; since $L$ is semisimple (the eigenvalues are distinct), it follows that $X_2$ must be diagonal. 

Then, choosing  $X_1=0$ and $D_1$ arbitrary we see that the diagonal part of $X_2$ must vanish as well. Thus the pairing is nondegenerate and the form $\o$ 
(\ref{defoGL}) is symplectic. \QED

The  corresponding Poisson structure is given by the following proposition.
\bp
\la{extPBKK}
The nonzero Poisson brackets corresponding to the symplectic form $\o$ are
\be
\la{PBGL}
 \{ G_{bj},G_{c\ell}\} =\frac { G_{b\ell} G_{cj} }{\l_j- \l_\ell}\ , 
\ \ j\neq \ell\;,\qquad \ \ 
 \{ G_{bk}, \l_\ell\} = -G_{bk} \delta_{\ell k}\;.
\ee
\ep
{\bf Proof.} 
The form \eqref{pariin} defines a map  $\Phi_{(G,L)} : T_{(G,L)}\Hcal \to  T_{(G,L)}^\star \Hcal$ given by 
\be
\label{phiGL}
\le\langle \Phi_{(G,L) }(X_1,D_1) , (X_2,   D_2) \ri\rangle := \o \big((X_1 ,D_1), (X_2, D_2)\big)
\ee
for all $(X_2,   D_2)$. Then (\ref{pariin}) implies
\be
\Phi_{(G,L) }(X,D)= \bigg( -D - [X,L],X^{^D}\bigg) \in T^\star_{(G,L)} \Hcal
\ee
where $X^{^D}$ and $X^{^{OD}}$ denote the diagonal and off-diagonal parts  of the matrix $X$, respectively  and the identification between a matrix and its dual is defined by  the trace pairing.  We denote  
\be
Q= -D - [X,L]\;,\hskip0.7cm \deltad=X^{^D}\;.
\la{QDD}\ee
Given now $(Q,\deltad)\in T^\star_{(G,L)}\Hcal$ we observe from   \eqref{QDD} that $D = -Q^{^D}$ and $X = \deltad + {\rm ad}^{-1}_L(Q^{^{OD}})$. The inverse of ${\rm ad}_{L}(\cdot) = [L,\cdot]$ is given explicitly by 
\be
{\rm ad}^{-1}_L(M)_{ab} = \frac {M_{ab}}{L_{aa} - L_{bb}}, \hskip0.7cm a\neq b
\ee
as a linear invertible map on the space of off--diagonal matrices. 

Thus $\Phi_{(G,L)}^{-1}: T_{(G,L)}^\star\Hcal \to  T_{(G,L)} \Hcal$ is given by
\be
\Phi_{(G,L)}^{-1} (Q,\deltad) = \bigg(\deltad + {\rm ad}_L^{-1} (Q^{^{OD}}), - Q^{^D}\bigg)
\ee
where $Q^{^{OD}}$ and $Q^{^D}$ denote the off-diagonal and diagonal parts, respectively.
The Poisson tensor $\mathbb P \in  \bigwedge^2 T_{(G,L)} \Hcal\simeq ( \bigwedge^2 T^\star_{(G,L)} \Hcal)^\vee$ is defined by  
\be
\la{PB}
\mathbb P_{(G,L)} \bigg((Q_1,\deltad_1), (Q_2, \deltad_2)\bigg):= 
\o \bigg(\Phi_{(G,L)}^{-1} (Q_1,\deltad_1),\Phi_{(G,L)}^{-1} (Q_2,\deltad_2)\bigg).
 \ee
Using the definition \eqref{pariin}, \eqref{phiGL} we get
$$
\mathbb P_{(G,L)}  \bigg((Q_1,\deltad_1), (Q_2, \deltad_2)\bigg)
$$
$$
=\tr \le(
-Q_1^{^D}\le(\deltad_2 + {\rm ad}_L^{-1} (Q_2^{^{OD}})\ri) 
+
Q_2^{^D}\le(\deltad_1 + {\rm ad}_L^{-1} (Q_1^{^{OD}})\ri)
 + {\rm ad}_L \le(\deltad_1 + {\rm ad}_L^{-1} (Q_1^{^{OD}})\ri)
 \le(\deltad_2 + {\rm ad}_L^{-1} (Q_2^{^{OD}})\ri)
  \ri) 
$$
which is equal to 
 $$
 \tr \Bigg(
Q_2^{^D}\deltad_1-Q_1^{^D}\deltad_2
 +Q_1^{^{OD}} {\rm ad}_L^{-1} (Q_2^{^{OD}})
  \Bigg) \;.
$$
To obtain the Poisson bracket between the matrix entries of $G$ and $L$ we now write
$Q = G^{-1} \d G$  and $\deltad = \d L = {\rm diag}(\d \l_1,\dots, \d \l_n)$.  

Choosing $Q_1 = \mathbb E_{jk}, \deltad_1 =0$ and $Q_2 = 0, \deltad_2 =\mathbb  E_{\ell\ell}$ we have
$$
(G^{-1})_{jb}\{ G_{bk}, \l_\ell\} =  \mathbb P( (G^{-1} \d G)_{jk} , \d \lambda_\ell) =- \delta_{jk} \delta_{\ell k}\  \Rightarrow \ \{ G_{bk}, \l_\ell\} =- G_{bk} \delta_{\ell k}\; .
$$
Choosing $Q_1 = \mathbb E_{ij}, Q_2 = \mathbb E_{k\ell},\  \ \deltad_1= \deltad_2 = 0$ we have
$$
\mathbb P_{(G,L)} \big((G^{-1} \d G)_{ij},(G^{-1} \d G)_{k\ell} \big)=
(G^{-1})_{ib} (G^{-1})_{kc}\{ G_{bj},G_{c\ell}\} = \frac {\delta_{jk} \delta_{i\ell}} { \l_j - \l_\ell}\;.
$$
\QED

\begin{proposition}
Introduce  the  $GL(n)$ {\it dynamical $r$-matrix}  (\cite{Etingof}, p.4):
$$
r(L)=\sum_{i<j} \f{E_{ij}\otimes E_{ji}-E_{ji}\otimes E_{ij}}{\l_i-\l_j}\;.
$$
where   $E_{ij}$ is an $n\times n$ matrix whose $(ij)$ entry equals $1$ while all other entries vanish.
Introduce also the matrix
$$
\Omega= \Omega_{\mathfrak {gl}(n)} :=\sum_{i=1}^n E_{ii}\otimes E_{ii}\; .
$$
Then the bracket (\ref{PBGL}) can be  written as follows:
\be
\{\m{G}^1,\m{G}^2\}=-\m{G}^1\m{G}^2 r(L)  \;,
\la{Dyno}
\ee
\be
\{\m{G}^1,\m{L}^2\}=-\m{G}^1 \Omega\;.
\la{Dyno2}
\ee

  \end{proposition}
 The proof is a straightforward computation.
 Notice that the formula (\ref{Dyno2}) can alternatively be written as follows:
\be
\{G,\l_j\}=-G E_{jj}\;.
\la{Dyno3}\ee
 
 The  Jacobi identity involving the brackets $\{\{\m{G}^1,\m{G}^2\},\m{G}^3\}$ implies (taking into account that $\ds \m{r}^{ij}=-\m{r}^{ji}$) the classical dynamical Yang-Baxter equation:
 (see (3) of \cite{Etingof}).
  \be
  [\m{r}^{12},\m{r}^{13}]+ [\m{r}^{12},\m{r}^{23}]+ [\m{r}^{23},\m{r}^{31}]+
  \sum_{i=1}^n\f{\ds \p \m{r}^{12}(L)}{\p \l_i} {{\m{E}^3}}_{ii}+\f{\ds \p \m{r}^{23}(L)}{\p \l_i} {\m{E}^1}_{ii}+
  \f{\ds\p \m{r}^{31}(L)}{\p \l_i} {\m{E}^2}_{ii}=0\;.
  \la{DYB}
  \ee
\begin{remark}\rm
 We did not find the  construction of this section in the existing literature. In the special case of the $SL(2)$ 
group, the Poisson algebra (\ref{Dyno}),  (\ref{Dyno2}) appeared in the work \cite{AlekFad} in the context of classical Poisson geometry of 
$T^*SL(2)$, see formulas (2),(3) in loc.cit.

As it was mentioned to us by L.Feher,  the Poisson structure (\ref{Dyno}), (\ref{Dyno2}) can be obtained from the canonical Poisson structure on $T^* SL(n)$ as follows. Consider an element $(G,A)\in T^* SL(n)$ and denote by $L$ the diagonal form of the matrix $A\in sl(n)$ (on an open part of the space where the matrix $A$ is diagonalizable).    The condition that $A$ is diagonal i.e $A=L$ is then a constraint of the second kind,
 according to Dirac's classification. The computation of the Dirac bracket for the pair $(G,L)$ starting from the canonical Poisson structure on  $T^* SL(n)$ leads to the Poisson structure (\ref{Dyno}), (\ref{Dyno2}), similarly to a computation given in \cite{Feher}. 
 
 \end{remark}

\subsubsection{Reduction to $SL(n)$}
To reduce to $SL(n)$ we observe that the proof of Prop. \ref{nondeg} holds also if we assume $\tr L=0$ and $\det G=1$. To compute the corresponding Poisson bracket we recall that inverting the restriction of a symplectic form to a symplectic submanifold is equivalent to the  computation of  the Dirac bracket. 

Let $h_1 :=\log  \det G$ and $h_2:=  \tr L$;  the Dirac bracket is then 
$$
\{F, H\}_D = \{F, H\} - \sum_{j=1}^2 \{F, h_j\} A_{jk} \{h_k, H\}
$$
where $A_{jk}$ is the inverse matrix to $\{h_j, h_k\}$: in our case we have 
$$
\{\log \det G, \tr L\} = -n\ \ \ \Rightarrow  \ \  A = \frac 1 n \le(\begin{array}{cc}
 0 &1\\
 -1 &0
\end{array}\ri)\;.
$$
Moreover a simple computation using \eqref{PBGL} shows that 
$$
\{G_{jk}, \det G\} =0\ ,\ \ \ \ \ \{G_{jk}, \tr L\} = -G_{jk} \;.
$$
Then  (we denote by $\{\}_{SL(n)}$ the Dirac bracket restricted to $\det G =1, \ \tr L =0$)
$$
\{G_{bj}, G_{c\ell}\}_{SL(n)}  = \{G_{bj}, G_{c\ell}\} \;,
$$
$$
\{\l_j, \l_k\}_{SL(n)}   = \{\l_j, \l_k\}=0\;,
$$
$$
\{G_{bk}, \l_\ell\}_{SL(n)}  = \{G_{bk}, \l_\ell\}  +  \frac 1 n \{G_{bk}, \tr L \} \{ \log \det G, \l_{\ell}\} 
 $$
 \be
 = - G_{bk} \delta _{\ell k}  + \frac 1 n G_{bk} =  G_{bk} \le(\frac 1n- \delta _{\ell k}\ri)\;.
 \label{SLnPB}
\ee
Equivalently the $SL(n)$  bracket is written as 
\be
\label{SPB}
\{\m{G}^1, \m{G}^2\}_{SL(n)} =  - \m{G}^1\m{G}^2 r(L)\;,
\qquad 
\{\m{G}^1, \m{L}^2\}_{SL(n)} = - \m{G}^1 \Omega
\ee
where now the matrix $\Omega$ is given by
$$
\Omega:= \Omega_{_{\mathfrak {sl }(n)}}= \sum_{j=1}^n E_{jj}\otimes E_{jj} - \frac 1 n \1 \otimes \1 = \sum_{j,k=1}^{n-1} (\mathbb A^{-1})_{jk}\;\alpha_j\otimes \alpha_k 
$$
and $\alpha_j = {\rm diag}(0,\dots, 1,-1,0,\dots)$ are the simple roots of $SL(n)$ and $\mathbb  A$ is the Cartan matrix of $SL(n)$; 
\be
\mathbb A = \le(
\begin{array}{cccccc}
2 & -1&0 \dots\\
-1&2&-1&0\dots\\ 
0 & -1 & 2&-1&\dots\\
0 & 0 &\ddots \\
\dots&&&-1 &2
\end{array}\ri)\;.
\la{Carmat}
\ee

\subsubsection{Relation to the Kirillov-Kostant bracket}

The Kirillov-Kostant bracket on $GL(n)$,  in tensor notation, takes the form 
\be
\{\m{A}^1,\m{A}^2\}=[\m{A}^1,P]=P(\m{A}^2-\m{A}^1)
\la{KK1}.
\ee
 Here $P$ is the permutation matrix of size $n^2\times n^2$ given by
\be
P=\sum_{i,j=1}^n E_{ij}\otimes E_{ij}\;.
\la{Perma}
\ee

 The regular symplectic leaves are the (co)adjoint orbits of diagonal matrices $L$ with distinct eigenvalues,  and on the orbit passing through $L$ the 
symplectic form of the Kirillov-Kostant bracket (\ref{KK1}) is equal to  (see \cite{BBT}, pp. 44, 45):
\be
\o_{KK}=-\tr \left(L G^{-1} \d G\wedge G^{-1} dG\right)
\la{wKK1}
\ee
where $G$ is any matrix diagonalizing $A$ i.e.  $A=GLG^{-1}$.
The form $\o_{KK}$ is  invariant under the transformation $G\to GD$ where $D$ is a diagonal matrix 
which may depend on $G$; such transformation leaves $A$ invariant. 
%\vskip2.0cm

\begin{theorem}
\la{thmrtoKK}
The map $(G,L)\mapsto A= GLG^{-1}$ is a Poisson morphism between the Poisson structure \eqref{PBGL} and  the Kirillov-Kostant Poisson  structure on $A$;
\be
\la{KKPB}
\{\tr (AF) ,\tr (AH) \}_{_{KK}}  = \tr  \bigg( A[H,F]\bigg)\;,\hskip0.7cm \forall F, H\in \mathfrak {gl}_n\simeq  \mathfrak {gl}_n^\vee
\ee
or, equivalently,
\be
\{\m{A}^1,\m{A}^2\}=[\m{A}^1,P]\;.
\la{KKalt}
\ee
\end{theorem}
\noindent
{\bf Proof.}
We have 
$$ \d A =  [\d G G^{-1} ,A] + G \d L G^{-1} = G\bigg([X,L] +\Lambda \bigg)G^{-1}\;.$$
Then 
$$
\{ (GLG^{-1})_{ab}, (GLG^{-1})_{cd}\} = \mathbb P_{(G,L)} \bigg(
\d (GLG^{-1})_{ab}, \d(GLG^{-1})_{cd}\bigg)
$$
\be
=\sum_{i,j,k,\ell=1}^{n}G_{ai}(G^{-1})_{jb} G_{ck} (G^{-1})_{\ell d}\mathbb P_{(G,L)} \bigg( \le([G^{-1}\d G,L] + \d L \ri)_{ij}, \le([G^{-1}\d G,L] + \d L \ri)_{k\ell}\bigg)\;.
\la{ers}
\ee
From the Poisson bracket (\ref{Dyno3}) we have
$$
\mathbb P( (G^{-1} \d G)_{jk} , \d \lambda_\ell) =- \delta_{jk} \delta_{\ell k}\\
\mathbb P \big((G^{-1} \d G)_{ij},(G^{-1} \d G)_{k\ell} \big)= \frac {\delta_{jk} \delta_{\ell i}} { \l_j - \l_\ell}\ ,\ \ j\neq \ell\;,
$$
\be
\mathbb P( \d\lambda_j , \d \lambda_\ell) =0\;.
\la{PGLGL}
\ee
Plugging \eqref{PGLGL} in \eqref{ers}  we see that the only terms giving non-trivial contributions   are the following:
$$
\sum_{i,j,k,\ell=1}^n
G_{ai}(G^{-1})_{jb} G_{ck} (G^{-1})_{\ell d}\mathbb P_{(G,L)} \bigg( \le([G^{-1}\d G,L] + \d L \ri)_{ij}, \le([G^{-1}\d G,L] + \d L \ri)_{k\ell}\bigg)
$$
$$
=\sum_{i,j,k,\ell=1}^n
G_{ai}(G^{-1})_{jb} G_{ck} (G^{-1})_{\ell d} (\l_i - \l_j)(\l_k -\l_\ell) \mathbb P \bigg( (G^{-1}\d G)_{ij},(G^{-1}\d G)_{k\ell}\bigg) 
$$
$$
=
\sum_{j,\ell=1}^n
G_{a\ell}(G^{-1})_{jb} G_{cj} (G^{-1})_{\ell d} (\l_\ell - \l_j) = A_{ad} \delta_{b c}  - A_{cb}\delta _{ad}\;.
$$

This expression coincides with  the Kirillov-Kostant Poisson bracket \eqref{KKPB}. 
\QED

A slight modification of this computation shows that the quadratic Poisson bracket (\ref{PBGL}) implies the Kirillov-Kostant bracket in $SL(n)$ case.

\subsection{Hamiltonian formulation of the  Schlesinger system in $G$--variables}
Consider the Schlesinger system written in terms of the matrices $G_j\in SL(n)$:
  \be
 \f{\p G_k}{\p t_j}= \f{A_j G_k}{t_j-t_k}\;,\;\;\;   j\neq k\;;\hskip0.7cm
 \f{\p G_k}{\p t_k}=-\sum_{k\neq j} \f{A_j G_k}{t_j-t_k}\;,\hskip0.7cm  \f{\p L_k}{\p t_j}=0
 \la{eqG}
 \ee
where
\be
A_j=G_j L_j G_j^{-1}\;.
\la{AjGj}
\ee
Matrices $L_j\in \mathfrak {sl}(n)$ are diagonal and the eigenvalues of $L_j$ are assumed to be distinct.

The  Poisson structure of the Schlesinger system (\ref{Schlsys}) written in terms of $A_j$ is known to be linear: it is based on  the  Kirillov-Kostant bracket.  
On the other hand, the hamiltonian formulation of the system (\ref{eqG}) involves the quadratic bracket defined by the dynamical $r$-matrix.

The following theorem can be  checked by direct calculation:
\begin{theorem}
\label{Sheh} Denote  by $\Acal_0 $ the space of matrices $\{(G_j,L_j)\}_{j=1}^N$ where $L_j$ are diagonal matrices with distinct eigenvalues.
Then the system (\ref{eqG}) is a multi-time hamiltonian system with respect to the Poisson structure 
on $\Acal_0 $
\be
\{\m{G_j}^1,\m{G_k}^2\}=-\m{G_j}^1\m{G_k}^2 r(L_k)\,\delta_{jk}\;,\hskip0.7cm
\{\m{G_j}^1,\m{L_k}^2\}=-\m{G_k}^1 \Omega\,\delta_{jk}\;.
\la{b33}
\ee
where $\delta_{jk}$ is the Kronecker delta.
The 
Hamiltonian defining the evolution with respect to "time" $t_k$ is given by
$$
H_k=\sum_{j\neq k} \frac{\tr A_j A_k}{t_k-t_j}\;,\hskip0.7cm k=1,\dots,N\;.
$$
\end{theorem}

We notice that for the Schlesinger system for matrices $A_j$ (\ref{Schlsys}) the Hamiltonians $H_j$ are the same as for the 
system (\ref{eqG}).

\subsection{Symplectic form and potential}

In the sequel we shall use the symplectic form associated to the bracket (\ref{b33}).
A direct computation     using the Poisson bracket \eqref{SLnPB} shows that the matrix $A= GLG^{-1}$ has the following Poisson brackets with $G$ and $L$:
\be
{\{A_{ab}, G_{jk}\} = G_{ak} \delta_{bj} - G_{jk} \frac {\delta_{ab}}n\;,}\hskip0.7cm  \{A_{ab}, \l_k\}=0\; .
\ee
Thus $\{\tr(X A), G\} = X G$ for any fixed matrix $X\in \mathfrak {sl}(n)$ and, therefore, the matrix $A = GLG^{-1}$ is  the moment map for the group action $G\mapsto S G$ on the space $\Hcal$. A similar statement, of course, holds for $GL(n)$ using the Poisson bracket \eqref{PBGL} instead.

Consider now  the diagonal  group action on $\Acal_0$ given by 
  \be
\{G_j,L_j\}_{j=1}^N\to \{S G_j,L_j\}_{j=1}^N
\la{Gracin}
\ee
where $S$ is an $SL(n)$ 
matrix. The previous computation shows immediately that the moment map corresponding to the group action $G_j\to SG_j$ on $\Acal_0$ is given by
\be
\{G_j,L_j\}_{j=1}^N\to {\mathfrak m} =\sum_{j=1}^N G_j L_j G_j^{-1}\;.
\la{momapin}
\ee

The space $\Acal$ is  defined by (\ref{defA}) as the space of the orbits  of the action (\ref{Gracin})
  of  in the zero level set of the moment map (\ref{momapin}).
This implies 
 the following theorem proven via  the standard symplectic reduction  \cite{Audin}:
\begin{theorem}\la{nondegmain}
The Poisson structure induced on $\Acal$ from the Poisson structure (\ref{b3}) on $\Acal_0$  via the  reduction 
on the level set $\sum_{j=1}^N G_j L_j G_j^{-1}=0$ of the moment map, corresponding to the group action $G_j\to SG_j$,
is non-degenerate and the corresponding symplectic form is given by 
\be
\o_{\Acal}=-\sum_{k=1}^N\tr (L_k G_k^{-1} \d G_k\wedge G_k^{-1} \d G_k) + \tr (\d L_k \wedge G_k^{-1} \d G_k)\;.
\la{wex2}
\ee
A symplectic potential $ \theta_\Acal$ for  $\o_{\Acal}$  is given by 
\be
  \theta_\Acal=\sum_{k=1}^N\tr (L_k G_k^{-1} dG_k)\;.
  \la{sympot2}
  \ee
 \end{theorem}

  \section{Monodromy symplectomorphism via Malgrange's form}
  \la{AMMalg}

We start from introducing the Malgrange form associated to a Riemann-Hilbert problem on an oriented graph and
discussing some of its properties, following \cite{Malg,Bert,Bert1}. From now on we work with the  $SL(n)$
case.

Let $\Sigma$ be an  embedded graph on $\C \mathbb P^1$ whose edges are  smooth oriented  arcs meeting transversally at the vertices. We denote by $\mathbf V$ the set of vertices of $\Sigma$. 
Consider  a  "jump matrix"  i.e. a function   $\Jum(z):\Sigma\setminus \mathbf V  \to SL(n)$ that  satisfies the following  properties
\begin{assumption}
\la{assumptionM}
\begin{enumerate}
\item In a small neighbourhood of  each point $z_0\in \Sigma \setminus \mathbf V $ the matrix $\Jum(z)$ is given by a germ of analytic function; 
\item for each $v\in \mathbf V $, denote by $\g_1,\dots, \g_{n_v}$ the edges incident at $v$ in a small disk centered thereof. Suppose first that all these edges are  oriented away from $v$ and enumerated in counter-clockwise order. 
Denote by $\Jum^{(v)}_j(z)$ the analytic restrictions of 
$\Jum$ to $\g_j$.
Assume that each $\Jum_j^{(v)}(z)$ admits an analytic extension to a full neighbourhood of $v$ and that these extensions satisfy the {\it local no-monodromy condition}
\be
\la{locnomonodromy}
\Jum_1^{(v)}(z)\cdots \Jum^{(v)}_{n_v}(z)=\1\;. 
\ee
If the  edge $\g_j$ is oriented towards $v$ then  $\Jum_j^{(v)}(z)$ is taken to be the inverse of   $\Jum(z)$. 
\end{enumerate}
\end{assumption}

Suppose now that the jump matrices form  an analytic family  depending on some deformation parameters and satisfying Assumption \ref{assumptionM}, and 
consider a  family of  Riemann-Hilbert problems on $\Sigma$.
\vskip0.3cm
{\bf Malgrange form for an arbitrary Riemann Hilbert Problem.} 
\label{RHPPhi}
 Let $\Phi(z):\C \mathbb P^1 \setminus \Sigma \to SL(n)$ be a  matrix--valued function, bounded everywhere and analytic on each face of $\Sigma$.  We also assume that the  boundary values on the two sides of each edge of $\Sigma$ are related by  
  \be
  \Phi_+(z)=\Phi_-(z) \Jum(z)\;, \ \ \ \forall  z\in \Sigma\setminus \mathbf V  \ ,\qquad \Phi(\infty)=\1\; .
  \la{RHprob}
  \ee
  where the ${+/-}$ boundary value is from the left/right, respectively, of the oriented edge. 
\vskip 3pt
  
\begin{definition} \cite{Malg}
The Malgrange 1-form on the deformation space of Riemann-Hilbert problems with given graph $\Sigma$ and jump matrices $\Jum$ is defined by 
\label{defMB}
\be
\label{omegaM}
\Theta [\Sigma,J]=\frac {1}{2i\pi} \int_{\Sigma} \tr \left(\Phi_-^{-1}\f{\d\Phi_-}{\d z} \d \Jum(z) \Jum^{-1}(z)\right) \d z
\ee
   where $\d \Jum$   denotes the total differential of $\Jum$ in the space of deformation parameters for  fixed $z$.
\end{definition}

  In (\cite{Bert1}, Thm. 2.1) it was proved the following formula for the exterior derivative of  (\ref{omMG}):
  \be
\la{dTheta}
\d \Theta [\Sigma,J]   = - \frac 1 2 \int_{\Sigma}  \ddz \tr\bigg(
\frac {\d}{\d z} (\d \Jum \Jum^{-1})  \wedge(\d \Jum \Jum^{-1})
\bigg)
+\Etavertex
\ee
where
 \be
\la{etavertex}
\Etavertex  =-\frac {1}{4i\pi}\sum_{v\in \mathbf V } \sum_{\ell=1}^{n_v-1} \tr \bigg( (\Jum_{\ell}^{(v)})^{-1}\d \Jum_{\ell}^{(v)} \wedge \d \Jum^{(v)}_{[\ell+1:n_v]} (\Jum^{(v)}_{[\ell+1:n_v]})^{-1} \bigg)\;.
\ee
Here the notation $\Jum^{(v)}_{[a:b]}$ stands for the product $J^{(v)}_{a} \cdots J^{(v)}_b$ for any two indices $a<b$.

In \cite{Bert1} the formula for $\Etavertex$ is written in a slightly different form and can be recast as the above expression by using the conditions \eqref{locnomonodromy}.

\paragraph{Malgrange form and Schlesinger systems.}   Let us now discuss how the form (\ref{omegaM}) can be used in the context of the Fuchsian equation (\ref{lsint}) and the associated Riemann-Hilbert problem.

Let $\mathbb D_j$ be small, pairwise non-intersecting disks centered at $t_j$, $j=1,\dots, N$. 

In order to define the inverse monodromy map  unambiguously, we need to fix the determination of the power in \eqref{asint}. To this end,  fix a point $\beta_j$ on the boundary of each of the disks $\mathbb D_j$  and declare that, within the disk $\mathbb D_j$, the power $(z-t_j)^{L_j}$ stands for $ |z-t_j|^{L_j} {\rm e}^{i\arg (z-t_j)L_j}$, where the argument is chosen between $\arg(\beta_j-t_j)$ and $\arg(\beta_j-t_j) + 2\pi$.  In particular the logarithm $\ln(z-t_j)$ is assumed to have the branch cut connecting $t_j$ with $\beta_j$ and the determination implied by the above. 

Choose now a collection of non-intersecting edges $l_1,\dots, l_N$ connecting $\infty$ with each of the $\beta_j$'s ($l_j$ is assumed to be transversal to  the boundary $\pa\mathbb D_j$ at $\beta_j$). 
Denote by $\Sigma$ the union of all the circles $\pa\mathbb D_j$ and the edges $l_j$.
Denote  by  $\mathscr D_\infty$ the  ``exterior'' domain which is  the complement of the union of the disks $\mathbb D_j$   and the graph $\Sigma$.

The solution  $\Psi(z)$ of \eqref{lsint} is a single--valued matrix function in $\scr D_\infty$ normalized by $\lim_{z\to\infty} \Psi(z) =\1$ where the direction lies within a sector lying between edges  $l_1$ and $l_N$.   Within each disk, with the above choice of determination of the logarithm, the analytic continuation of $\Psi$ has the local expression \eqref{asint}. The  ``connection matrices'' $C_j$ are uniquely determined by a choice of $G_j$ and the determination of the logarithm. 

We have therefore defined the (extended) {\it monodromy map}
\be
\label{FT}
 \le\{(G_j, L_j):\ \ \ \sum_{j=1}^{N} G_jL_j G_j^{-1} =0\ri\} \to \le \{(C_j, L_j\}: \ \ \prod_{j=1}^N C_j {\rm e}^{2i\pi L_j} C_{j}^{-1} = \1\ri\}.
\ee
Although this monodromy map depends on $\Sigma$ and the determinations of the logarithms, we are not going to indicate it explicitly.

 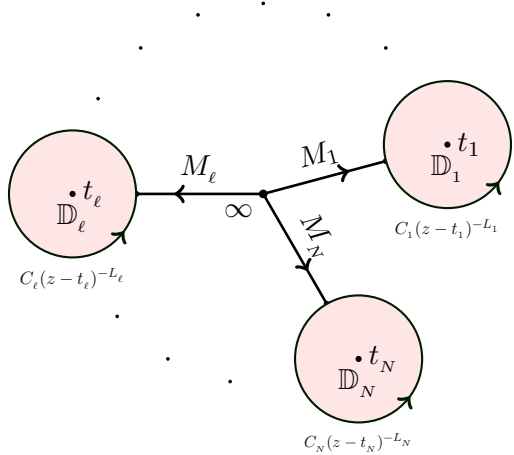
\begin{figure}[htb]
 \begin{center}
 \begin{tikzpicture}[scale=1.69]
\draw[line width = 1pt, postaction={decorate,decoration={{markings,mark=at position 0.7 with {\arrow[black,line width=1.5pt]{>}}}} }]
(0,0)--node[above, sloped]{{$ M_1$}}(15:1);

\draw[line width = 1pt, postaction={decorate,decoration={{markings,mark=at position 0.7 with {\arrow[black,line width=1.5pt]{>}}}} }]
(0,0)--node[above]{{$ M_{_\ell}$}}(180:1);

\draw[line width = 1pt, postaction={decorate,decoration={{markings,mark=at position 0.7 with {\arrow[black,line width=1.5pt]{>}}}} }]
(0,0)--node[above,sloped ]{{$ M_{_N}$}}(-60:1);
\draw [fill] circle [radius=0.03];
\node [below left] at (0,0) {$\infty$};
\foreach \x in {50,70,90,110,130,150, 220, 240,260}
\draw [fill] (\x:1.5)  circle [radius=0.01];

\draw [fill] (15:1)  circle[radius=0.02];
\draw [fill] (-60:1)  circle[radius=0.02];
\draw [fill] (180:1)  circle[radius=0.02];

\begin{scope}[shift={(15:1.5)}]
\draw[green!10!black,  fill = red!10!white, thick,postaction={decorate,decoration={markings,mark=at position 0.9 with {\arrow[line width=1.5pt]{>}}}}] (0:0.5) arc (0:360:0.5) node [pos = 0.75, below,sloped] {\resizebox{1.4cm}{!}{$C_{_1} (z-t_{_1})^{-L_{_1}}$}};

\node [right] at (0:0)  {$t_1$} ;
\draw [fill](0,0)  circle[radius=0.02];
\node [below] at (0:0)  {$\mathbb D_{_1}$} ;

\end{scope}

\begin{scope}[shift={(180:1.5)}]
\draw[green!10!black,  fill = red!10!white, thick,postaction={decorate,decoration={markings,mark=at position 0.9 with {\arrow[line width=1.5pt]{>}}}}] (0:0.5) arc (0:360:0.5) node [pos = 0.75, below,sloped] {\resizebox{1.4cm}{!}{$C_{_\ell} (z-t_{_\ell})^{-L_{_\ell}}$}};
\node [below] at (0:0)  {$\mathbb D_{_\ell}$} ;

\node [right] at (0:0)  {$t_{_\ell}$} ;
\draw [fill](0,0)  circle[radius=0.02];
\end{scope}

\begin{scope}[shift={(-60:1.5)}]
\draw[green!10!black, fill = red!10!white, thick,postaction={decorate,decoration={markings,mark=at position 0.9 with {\arrow[line width=1.5pt]{>}}}}] (0:0.5) arc (0:360:0.5) node [pos = 0.75, below,sloped] {\resizebox{1.45cm}{!}{$C_{_N} (z-t_{_N})^{-L_{_N}}$}};

\node [right] at (0:0)  {$t_{_N}$} ;
\node [below] at (0:0)  {$\mathbb D_{_N}$} ;
\draw [fill](0,0)  circle[radius=0.02];
\end{scope}

\end{tikzpicture}
\end{center}
\caption{Graph $\Sigma$ and jump matrices on its edges used in the calculation of the form $\Theta$. }
\label{FigSchles}
\end{figure}

An example of  the graph  $\Sigma$ is  shown  in Fig.\ref{FigSchles}; the graph looks like $N$ "cherries" whose "stems" are attached to the point $z=\infty$. Introduce the piecewise analytic matrix on its faces as follows
\be
\la{PhiPsi}
\Phi(z) = \left\{
\begin{array}{cc}
\Psi(z)\;, & z\in \mathbb D=\C \mathbb P^1\setminus \Sigma \setminus\bigcup_{j=1}^N \mathbb D_j\;; \\[5pt]
\Phi_j(z) := \Psi(z) C_j (z-t_j)^{-L_j}  \;, & z\in \mathbb D_j\; .
\end{array}
\right.
\ee
The function  $\Phi$ solves a Riemann--Hilbert Problem on $\Sigma$ with the   jump matrices on its edges  indicated in Fig.\ref{FigSchles}:

\be
\la{JAMA}
J = \left\{
\begin{array}{cc}
M_j = C_j {\rm e}^{2i\pi L_j} C_j^{-1}\;,& z\in l_j \;; \\[5pt]
 C_j (z-t_j)^{-L_j} \;, & z\in \pa \mathbb D_j 
\end{array}
\right.
\ee
where $l_j$ is the "stem" of the $j$th cherry.

The matrix function $\Phi$ given by (\ref{PhiPsi})  is the unique solution of the Riemann--Hilbert problem with jump matrices (\ref{JAMA}):
\be
\label{RHPphi}
\Phi_+(z) = \Phi_-(z) J(z)\; , \hskip0.7cm 
\Phi(\infty)=\1,
\ee

The solution  of the Riemann-Hilbert problem   exists for generic set of data $\{C_j,L_j,t_j\}$; this
 solution  provides the inverse of the map  in \eqref{FT}. 
 We  emphasize that the inverse monodromy map depends on the isotopy class of $\Sigma$ and on the fixing of the branches of the  logarithms. 

In the context of Fuchsian systems the general Malgrange form  in Def. \ref{defMB} specializes to the following definition:
\begin{definition}
\label{defThetaM}
The Malgrange one form $\Theta \in T^*_p \widetilde \Mcal$ is the  form  defined by the expression \eqref{omegaM}, where $\Phi$ is the solution of the Riemann--Hilbert problem  \eqref{JAMA}, \eqref{RHPphi}.  
\end{definition}
It is known \cite{Malg} that the form $\Theta$ is a meromorphic form on $\widetilde \Mcal$; the set of poles 
of $\Theta$ is called the "Malgrange divisor"; on this divisor the 
 Riemann-Hilbert problem fails to have a solution. Moreover, the residue along this divisor is a {\it positive} integer \cite{Malg}.

The deformation parameters involved in the expression \eqref{omegaM} for $\Theta$ are $C_j, L_j$ subject to the monodromy relation  $\prod_{j=1}^N C_j {\rm e}^{2i\pi L_j} C_j^{-1}=\1$, and the locations of the poles $t_1,\dots, t_N$. 

\begin{theorem}\la{ththa}
The  form $\Theta\in T^*\widetilde \Mcal$  (\ref{omegaM}) and the potential $\widetilde \theta_{\Acal}\in T^*\widetilde \Acal$   are related by 
\be
\label{omMG}
 \Theta   
 =  
(\widetilde{ \Fcal}^{-1})^* \le(\theta_\Acal - \sum_{j=1}^N H_j \d t_j \ri)
\ee
where $H_j$ are the Hamiltonians (\ref{Hkint}).  
Denote now by $\pa_{t_j} \in T\widetilde \Mcal$ the vector field of differentiation w.r.t. $t_j$ keeping the monodromy data constant. Then the contraction of $\Theta$ with  $\pa_{t_j}$  is given by
\be
\label{Thetat}
\Theta (\pa_{t_j}) =  H_j\,.
\ee
Equivalently, the contraction of (the pullback via the inverse monodromy map of) $\widetilde \theta_\Acal$ with $\pa_{t_j}$ equals $2H_j$. 
\end{theorem}

{\it Proof.} 
The simplest way to prove (\ref{omMG}) is via the {\it localization formula} \cite{ILP}  using the
  Riemann-Hilbert problem  defined on the graph $\Sigma$ shown in Fig.\ref{FigSchles}.
  To simplify the notation we will not indicate explicitly the pullbacks, but simply consider the matrices $G_j$ as functions of times and monodromy data via the inverse monodromy map.

In the formula (\ref{omegaM}) the function $\G_-$ coincides with the boundary value of  the solution, $\Psi$, of the ODE \eqref{lsint} in  the domain $\mathbb D$. Therefore, denoting $\d\Phi/\d z$ by $\Phi'$ we have:
\be
\la{36}
\tr \le(\Phi_-^{-1} \Phi_-' \d \Jum\Jum^{-1} \ri) = 
\tr \le( A(z) \Phi_-\d \Jum\Jum^{-1}\Phi_-^{-1} \ri)  \;.
\ee
Here we have  used the fact that $\G_-$ coincides with $\Psi$ and therefore  $\Phi_-'\Phi_-^{-1} = A(z)$. 
Moreover  we have
$$
\Phi_-\d \Jum\Jum^{-1}\Phi_-^{-1} = \d\le(\Phi_- \Jum\ri) \Jum^{-1}\Phi_-^{-1}  - \d \Phi_- \Phi_-^{-1} = \d \Phi_+ \Phi_+^{-1} -\d \Phi_- \Phi_-^{-1} 
$$
since $\G_+ = \G_- J$.
Thus  \eqref{omegaM} can be equivalently written as follows
\be
\Theta= \f{1}{2\pi i} \int_\Sigma \tr \le( A(z)(\d \Phi_+ \Phi_+^{-1} -\d \Phi_- \Phi_-^{-1} )\ri) \d z\;
\label{oooo}
\ee
and further represented as
\be\la{thsplit}
\Theta= \f{1}{2\pi i} \int_{\p \mathbb D} \tr (A(z) \d \Psi \Psi_-^{-1}) \d z+ \f{1}{2\pi i} \sum_{j=1}^N \int_{\p \mathbb D_j}
\tr (A(z)\d \Phi_+ \Phi_+^{-1})\d z\; .
\ee
The first integral in the r.h.s. of (\ref{thsplit}) vanishes since the integrand is holomorphic in $\mathbb D$.
Thus  \eqref{thsplit} reduces to (this is the expression that also appears in \cite{ILP}, formula (1.11)):
\be
\Theta= \sum_{j} \res{z=t_j}  \tr \le(A(z) \d \G_j(z) \G_j^{-1}(z) \ri)\d z\,.
\la{thetaint}
\ee

The expression (\ref{thetaint}) can be further evaluated in the coordinate system given by $(C_j,L_j,t_j)$. Namely, the contribution of derivatives with respect to monodromy data
$(C_j,L_j)$ into (\ref{thetaint}) is obtained by evaluation of  $\d \G_j(z) \G_j^{-1}(z)$  at the poles $t_j$ which gives the monodromy part of $\tt_{\Acal}$ in (\ref{omMG}). 

A straightforward local analysis using \eqref{PhiPsi} shows that:
$$
\pa_{t_k} \G_j(z) \G_j^{-1}(z)\bigg|_{z=t_j}  = \pa_{t_k} G_j G_j^{-1}  - \delta_{kj}  \pa_{t_k} G_kG_k^{-1}  - \delta_{kj}  [A_k, \G_j'(t_j)\G_j(t_j)^{-1}]\;.
$$
Thus 
$$
\Theta=  \sum_{j} \res{z=t_j} \tr \bigg( A(z) \d \G_j(z)\G_j^{-1} (z) \bigg)  dz
= \sum_{j} \tr \bigg( A_j \d G_jG_j^{-1} \bigg) 
- \sum_{j} \d t_j \tr \bigg(A_j \pa_{t_j} G_j G_j^{-1}    \bigg)\;.
$$
Finally, due to the  Schlesinger equations for $G_j$ \eqref{eqGint}
we get 
$$
\Theta= \sum_{j} \tr \left( A_j \d G_jG_j^{-1} \right) 
- \sum_{j} \d t_j \sum_{ k\neq j}\frac{\tr A_j A_k  }{t_j-t_k}\;.
$$
Recalling that the Jimbo-Miwa Hamiltonians are given by  $H_j=  \sum_{ k\neq j}\frac{\tr A_j A_k   }{t_j-t_k}$ and that the first term equals the potential $ \widetilde {\theta}_\Acal$  on $\At$, we
arrive at (\ref{omMG}).

As a corollary of the Schlesinger equations (\ref{eqGint})  the contraction of $\widetilde  \theta_\Acal$ with a vector field $\pa_{t_j}$  (for fixed monodromy data) is 
$$
\widetilde  \theta_\Acal (\pa_{t_j}) = 2 H_j\;.
$$
Therefore, the total $\d  t_j$ -  part of the form $\Theta$ for fixed monodromies equals to $\sum_{j=1}^N H_j\d t_j$.
\QED

\paragraph{\bf Symplectic form on the  monodromy manifold.}
 We start from defining the two-form on the monodromy manifold which is one of central objects of this paper.

\begin{definition}
Define  the following 2-form on $\Mcal$ (\ref{defM}):

  \be
   \o_\Mcal=\f{1}{4\pi i}(\o_1+\o_2)
   \la{sympmon}
   \ee
   where
   \be
\o_1= \sum_{\ell=1}^{N} 
\tr \left(  M_\ell^{-1} \d M _{\ell}  \wedge K_{\ell }^{-1} \d K_{\ell}\right) 
+  \sum_{\ell=1}^N  \tr\left( \Lambda^{-1}_\ell C_{\ell} ^{-1}\d C_{\ell}\wedge  \Lambda_\ell  C_\ell ^{-1}\d C_\ell \right)\;,
\la{w1main}
\ee
\be
\o_2=2\sum_{\ell=1}^N  \tr  \left(\Lambda_\ell^{-1} \d \Lambda_\ell \wedge C_\ell^{-1} \d C_\ell
\right)
\la{w2main}\ee
and $K_{\ell}=M_1\dots M_\ell$.
\end{definition}

 On the monodromy manifold $M_1\dots M_N=\1$ the form $\o_\Mcal$ is invariant under simultaneous transformation
$C_j\to S C_j$ with $S$ is an arbitrary $SL(n)$-valued function on $\Mcal$.

\br
The restriction of the form $-2i\pi \o_\M$  on the leaves $\Lambda_j = $ constant (under such restriction  $\o_2=0$ and hence $-2i\pi\o_\Mcal = -\o_1/2$) coincides with  the symplectic
form  on the symplectic leaves of  the $GL(n)$ Goldman bracket found in  \cite{AlekMal} (formula 
(3.14); the  case of this formula relevant for us corresponds to $k=2\pi$ and  $g=0$ in the notation of \cite{AlekMal}). 

As we prove below in Corollary \ref{cornond}, the form $\o_\Mcal$ is non-degenerate on the space $\Mcal$, which is
a torus fibration (with fiber the product of $N$ copies of the $SL(n)$ torus of diagonal matrices)  over the union of all the symplectic leaves of the Goldman bracket. The fact that $\M$ is a torus fibration  is simply due to the fact that the fibers of the map $(C_j,\Lambda_j) \to M_j = C_j \Lambda_j C_j^{-1}$ are obtained by multiplication of the $C_j$'s  on the right by diagonal matrices.
\er
   
Let us trivially extend the form $\o_\Mcal$ to the space $\widetilde {\Mcal}$ (\ref{Mt}) which includes also the variables $t_j$.
 This extension is denoted by $\wt_\Mcal$.

 {\bf Relation between forms $\Theta$ and $\omega_{\Mcal}$.}
   The following theorem was stated in  \cite{Bert} in slightly different notations without direct proof. The  proof is given below.
\begin{theorem}\la{dTTH} 
The exterior derivative of the  form $\Theta$ is given by the pullback of the form
$\widetilde \o_\Mcal$  (\ref{sympmon}) under the monodromy map:
  \be
  \d\Theta=[\widetilde {\Fcal}]^*\widetilde {\o}_\Mcal\;.
  \la{Golmalg}
  \ee
  \end{theorem}
 
 {\it  Proof.}
   Let us  apply the formulas (\ref{dTheta}), (\ref{etavertex}) to the graph $\Sigma$  depicted in Fig.\ref{FigSchles} 
  with  indicated jump matrices.
  The integral over $\Sigma$ in the formula \eqref{dTheta} then reduces to a sum of integrals over $\pa \mathbb D_\ell$'s because the jump matrix $J(z)$ on the cuts is constant with respect to $z$. We denote by $\beta_\ell$ the 
three-valent vertices where the circles around $t_\ell$ meet with the edges going towards $z_0$. 
Let us consider the contribution of one of  the integrals over $\pa \mathbb D_\ell$ to (\ref{dTheta}).

We will drop the index $_\ell$ for brevity in the formulas below. Notice also that  $\d L \wedge \d L=0$ because the matrix $L$ is diagonal. Letting $J(z)  =C (z-t)^{-L}$  we get
\be
-\frac 1 2 \oint \ddz \tr\le(\frac {\d}{\d z}\Big( \d J(z) J(z)^{-1} \Big)\wedge\d  J(z) J(z)^{-1}  \ri)
=
\frac 1 2 \le(
\frac {\d t\wedge L\d L}{(\beta-t)}
+ 
\d L \wedge C^{-1}\d C  
\ri)\;.
\label{betaeta2}
\ee
In the course of the computation  we have used that 
$$
\int_{\beta}^\beta \ddz \frac {\log (z-t)}{(z-t)^2} = -\frac 1 {\beta-t}
$$
where the integration goes along the circle $|z-t|=|\beta-t|$ starting at $z=\beta$.
We now turn to the evaluation of the term $\Etavertex$  \eqref{etavertex}. The set of  vertices $\mathbf V $ consists of $\mathbf V = \{z_0, \beta_1,\dots, \beta_N\}$. The contribution coming from the vertex $z_0$ is precisely the first term in $\omega_1$ \eqref{w1main} (in \eqref{w1main} this term is  simplified using the local no-monodromy condition \eqref{locnomonodromy}).

To evaluate the contribution of the vertex $\beta=\beta_\ell\in \mathbf V $ we observe that this vertex is tri-valent and the jump matrices on the  three incident arcs are 
$$
J_1 = C {\Lambda^{-1}}  C^{-1} ,\qquad J_2 = C (\beta-t)^{-L} \ ,\qquad J_3 = (\beta-t)^L{\rm e}^{2i\pi L}  C^{-1}
$$
where $\Lambda:= {\rm e}^{2i\pi L}$.
In the definition it is assumed that $(z-t)^L$ is defined with a branch cut extending from $t$ to $\beta$. 
Since $J_1J_2J_3 =\1$ the contribution of the vertex  to \eqref{etavertex} reduces to  the term 
$$\frac {-1}{4i\pi} \tr\le(J_1\d J_2\wedge \d  J_3\ri) =\frac {-1}{4i\pi} \tr\le(J_2^{-1}\d J_2\wedge \d  J_3J_3^{-1}\ri) \;. $$
Recall that $L, \Lambda$ are diagonal; we have then
$$
J_2^{-1} \d J_2 =(\beta-t)^L C^{-1}  \d C (\beta - t)^{-L}  + (\beta-t)^L \frac {L \d t}{\beta-t} (\beta- t)^{-L}  -  \d L \log (\beta-t)\;,
$$
\be
\d J_3J_3^{-1}  = 
\frac{-\d t L}{(\beta - t)}  
 + \le(\log(\beta-t) + 2i\pi \ri) \d L
   - (\beta-t)^L \Lambda C^{-1} \d C \Lambda^{-1} (\beta-t)^{-L}
   .\la{dJs}
\ee

Then a straightforward computation gives 
$$
\frac {-1}{4i\pi}\tr\le(J_2^{-1}\d J_2\wedge \d  J_3J_3^{-1}\ri) 
$$
 \be
=
\frac {-1}{4i\pi}\tr \bigg(
C^{-1}  \d C \wedge  \Lambda^{-1}\d \Lambda
-
C^{-1}  \d C \wedge  \Lambda C^{-1} \d C \Lambda^{-1}
+
2i\pi \frac {L \d t}{\beta-t}  \wedge  \d L
 \bigg) \;.%%%%%%%%%%
 \la{dJss}
\ee
Summing up \eqref{betaeta2} (the contribution of the integral) with \eqref{dJss} (the contribution coming from the vertex $\beta=\beta_\ell$) 
we get
 $$
 \eqref{betaeta2} + \eqref{dJss}= 
\frac {1}{4i\pi}\tr \bigg(
-2C^{-1}  \d C \wedge  \Lambda^{-1}\d \Lambda
+
C^{-1}  \d C \wedge  \Lambda C^{-1} \d C \Lambda^{-1}
 \bigg)\;.
$$

Then summing over all contributions from  vertices $\beta_\ell$ leads to (\ref{sympmon}).

Summarizing, the first term in (\ref{w1main}) corresponds to the $N$-valent vertex.  The second term in  (\ref{w1main}) 
together with the term (\ref{w2main}) arise from the contributions of cherries and the three-valent vertices 
formed by cherries and their stems.  \QED
  
  \vskip 5pt

  This theorem immediately implies the following corollary, which can also be deduced from previous results of \cite{Boalch2}.  
  \begin{corollary}
  \la{cornond}
  The form $\o_\Mcal$ (\ref{sympmon}) is closed and non-degenerate on the  monodromy 
  manifold $\Mcal$.
  \end{corollary}
  
  \paragraph{ Strong version of Its-Lisovyy-Prokhorov conjecture.}
  
  The theorem \ref{dTTH} proves the  "strong" version of the  ILP conjecture  (\ref{ILP}).
  To state this conjecture in the present setting we consider the form (1.11) or (2.7) of \cite{ILP} which 
  we denote by $\Theta_{ILP}$ to avoid  confusion with the  notations of this paper (see also the identity (\ref{tpart}) below):
 \be
  \Theta_{ILP}=\sum_{j<k}^N\tr A_j A_k d\log (t_j-t_k)+\sum_{j=1}^N \tr (L_j G_j^{-1} \d_{\Mcal} G_j)\;.
  \la{strongomega}
  \ee
The Conjecture from section 1.6 of \cite{ILP} refers to the restriction of the form to the symplectic leaves $L_j=$constants. We refer to this as the {\it weak Its-Lisovyy-Prokhorov conjecture}; in  this formulation $\d_{\Mcal}$ refers to the differential only with respect to the connection matrices $C_j$.    This  "weak" version of the conjecture is proved on the basis of known results  \cite{Hitchin,AlekMal2,KorSam}  in the next section.

  The statement of Theorem \ref{dTTH}  is the {\it strong} version of the above conjecture: in this version the differential $\d_{\Mcal}$ is with respect to all monodromy data including the $L_j$'s.

\paragraph{Generating function of the monodromy map.} 
  The closure of $\omega_\Mcal$ guarantees the  local existence of a {\it symplectic potential}. 
  Denoting any  such local potential by $\theta_\Mcal$ (such that $\d\theta_\Mcal=\o_\Mcal$)  we define  the (local on $\Mt$) generating function $\Gcal$ 
  as follows
  \be
  \d\Gcal=\sum_{k=1}^N\tr (L_k G_k^{-1} \d G_k)-\sum_{j=1}^N H_k \d t_k - \widetilde \theta_\Mcal
  \la{defcon}
  \ee
  where $G_k$ and $H_k$ are considered as functions on $\Mt$ under the inverse monodromy map.
  
The equation (\ref{defcon}) can be used to extend the definition of Jimbo-Miwa tau-function to include its dependence on monodromies. 
  Irrespectively of the choice of $\theta_\Mcal$,  the formula \eqref{Thetat}  implies the following theorem
 \begin{theorem}\la{tauGmain}
For any choice of symplectic potential $\theta_\Mcal$ on $\Mcal$ the dependence of the generating function  $\Gcal$ (\ref{defgen}) on $\{t_j\}_{j=1}^N$ coincides with $t_j$-dependence of the 
isomonodromic Jimbo-Miwa tau-function. In other words, $e^{-\Gcal}\tau_{JM}$ 
depends only on monodromy data $\{C_j,L_j\}_{j=1}^N$. 
\end{theorem}
          
    In Section \ref{tau_fun_mon} we are going to use this theorem to {\it define }
    the isomonodromic tau function as exponent of the generating function $G$    
    under a special choice of the symplectic potential $\theta_\Mcal$ based on the use of Fock-Goncharov coordinates.

\br \rm
 "Extended" character varieties with non-degenerate symplectic form were considered  in the '94 paper \cite{Jeffrey}  and later in the  paper \cite{Boalch2}.
 In  (\cite{Boalch2} Corollary 1) it was proven that the pullback of a symplectic form from the extended monodromy manifold coincides with a symplectic form on $(L_j, G_j)$ side. 
 The description of the corresponding Poisson bracket, construction of 
symplectic potentials, Malgrange form, the tau-function and coordinatization in term of Fock-Goncharov parameters were not considered before,
to the best of our knowledge.
\er

 \section{Standard monodromy map and weak version of Its-Lisovyy-Prokhorov conjecture}
\la{previous}
  
  Here we show that a weak version of Its-Lisovyy-Prokhorov conjecture 
  can be derived in a simple way from previous results of \cite{Hitchin,AlekMal2} or \cite{KorSam} where a symplectomorphism 
  between the space of coefficients $\{A_j\}$ with given set of eigenvalues of the Fuchsian equation 
  (\ref{lsint}) and a symplectic leaf of Goldman bracket was proved.
  
  First, consider the  submanifold $\Acal_L$ of $\Acal$ such that  the diagonal form of each of the matrices $A_j$
  is fixed:
  \be
\Acal_L=\le \{\{A_i\}_{i=1}^N,\;\; A_i\in {\mathcal O}(L_i)\;,\;\;\sum _{i=1}^N A_i=0\ri\}/\sim
\la{A1}
\ee
where $\sim$ is the equivalence over simultaneous adjoint transformation $A_i\to SA_iS^{-1}$ of all $A_i$ for $S\in SL(n)$; $L=(L_1,\dots,L_N)$ where $L_j$ is the diagonal form of $A_j$ and ${\mathcal O}(L)$ is the (co)-adjoint orbit of the diagonal matrix $L$.  We assume that diagonal entries of each $L_j$ do not differ by an integer.
  
  Consider similarly also the space $\Mcal_L$ which is the subspace of the $SL(n)$ character variety of $\pi_1(\CP1\setminus\{t_j\}_{j=1}^N)$ such that the diagonal form of the matrix $M_j$ equals to $\Lambda_j=e^{2\pi i L_j}$.

 The Kirillov-Kostant brackets (\ref{KK1}) for each $A_j$:
 \be
 \{{\m{A}^1}_j,{\m{A}^2}_k\}=[{\m{A}^1}_j,P]\;\delta_{jk}
 \la{KK3}
 \ee
  can be equivalently rewritten in the $r$-matrix form
  \be
\{\m{A}^1(z){\,,\, }\m{A}^2(w)\} = 
\frac1{z-w}\,[P, \m{A}^1(z) + \m{A}^2(w)]\;.
\la{Rmat}
\ee
The Schlesinger equations for $A_j=G_j L_j G_j^{-1}$ which follow from the system (\ref{eqGint}) for $G_j$ take the form:
\be
 \f{\p A_k}{\p t_j}= \f{[A_k, A_j]}{t_k-t_j}\;,\hskip0.7cm   j\neq k\;; \hskip0.7cm
 \f{\p A_j}{\p t_j}=-\sum_{k\neq j} \f{[A_k, A_j]}{t_k-t_j}\;.
 \la{Schles}
 \ee
These equations are Hamiltonian,
$$
\f{\p A_k}{\p t_j}=\{H_j, A_k\}\;,
$$
with the  Poisson structure  given by (\ref{Rmat}) and the (time dependent) Hamiltonians $H_j$  defined by (\ref{Hkint}). Notice that these Hamiltonians commute  $\{H_k,H_j\}=0$ and satisfy the equations $\pa_{t_k} H_j = \pa_{t_j} H_k$. 

After the symplectic reduction to the space of orbits of the global $Ad_{GL(N)}$ 
  action and restriction to the level set $\sum_{j=1}^N A_j=0$ of the corresponding 
  moment map one gets a degenerate  Poisson structure; its symplectic leaves coincide with $\Acal_L$ \cite{Hitchin}.
  The symplectic form on $\Acal_L$ can be written as
  \be
  \o^L_{\Acal}=-\sum_{k=1}^N\tr (L_k G_k^{-1} dG_k\wedge G_k^{-1} dG_k) \;.
  \la{wLA}
    \ee
  The form (\ref{wLA}) is independent of the choice of matrices $G_j$ which diagonalize $A_j$; moreover, it is invariant under simultaneous transformation $A_j\to S A_j S^{-1}$ and  thus it is indeed defined on the space $\Acal_L$. 
  
  The  $SL(n)$ character variety is  equipped with the  Poisson structure given by the Goldman bracket  defined as follows (see p.266 of
\cite{Goldman1}):  for 
  any two loops $\sigma,\widetilde {\sigma}\in \pi_1(\CP1\setminus\{t_i\}_{i=1}^N)$ the Poisson bracket between the traces of the corresponding monodromies  is given by
  \be
\Big\{\tr M_\sigma,\;\tr M_{\widetilde{\sigma}}\Big\}_G= \sum_{p\in \sigma\cap \widetilde{\sigma}} \nu(p)\,\left( \tr (M_{\sigma_p\widetilde{\sigma}})
-\f{1}{n}  \tr M_\sigma  \tr M_{\widetilde{\sigma}}  \right)\;.
\la{Goldmanbr}
\ee
where $\nu(p)=\pm 1$ is the contribution of point $p$ to the intersection index of $\sigma$ and $\widetilde{\sigma}$.

The space $\Mcal_L$ is a symplectic leaf of the $SL(n)$ Goldman bracket;  
the Goldman's symplectic form on $\Mcal_L$ coincides with   $-\f{1}{2}\o_1$ \cite{AlekMal} where $\o_1$ is defined in  (\ref{w1main}). We define
 \be
  \o_\Mcal^L=\f{1}{4\pi i} \o_1\;.
  \la{wLG}
  \ee
The study of the symplectic properties of the map (\ref{Fdef}) was initiated in \cite{Hitchin,AlekMal2,KorSam}. In \cite{Hitchin,AlekMal2} two  different proofs  were given
of the fact  that the monodromy map $\Fcal^t$ is a symplectomorphism   i.e.
\be
(\Fcal^t)^* \o_\Mcal^L =  \o_\Acal^L\;.
\la{HitAM}
\ee

  In \cite{KorSam} the brackets between the monodromy matrices themselves were obtained starting from (\ref{Rmat}); the result is given by 
  \be
\{{\m{M}^1}_i, {\m{M}^2}_i \}^*=\pi i\,P( {\m{M}^1}_i{\m{M}^1}_i- {\m{M}^2}_i{\m{M}^2}_i)\;,
\la{KS1}
\ee
\be
\{{\m{M}^1}_i , {\m{M}^2}_j\}^* =
\pi i \,P \,\Big( {\m{M}^1}_j  {\m{M}^1}_i + {\m{M}^2}_i  {\m{M}^2}_j    - 
 {\m{M}^1}_i  {\m{M}^2}_j - {\m{M}^1}_j  {\m{M}^2}_i \Big)\;, \hskip0.7cm i< j 
\label{KS2} \ee
where $P$ is the matrix of permutation of two spaces.
The brackets  (\ref{KS1}), (\ref{KS2}) were computed for the basepoint $z_0=\infty$
on the level set $\sum_{j=1}^N A_j=0$ of the moment map; thus the algebra 
(\ref{KS1}), (\ref{KS2}) does not satisfy the Jacobi identity. However, the Jacobi identity is restored for the algebra of $Ad$-invariant objects i.e. for traces of monodromies; moreover, for any two loops $\sigma$ and $\widetilde{\sigma}$ we have
(\cite{Schemmel}; see also Thm. 5.2 of \cite{Marta1} where this statement was proved
for $n=4$, $N=2$ case): 
\be
\{\tr M_\sigma,\tr M_{\widetilde{\sigma}}\}^*=-2\pi i\{\tr M_\sigma,\tr M_{\widetilde{\sigma}}\}_G
\la{Goldman1}
\ee
  which gives an alternative proof of (\ref{HitAM}).

  Let us now show that (\ref{HitAM}) implies the weak version of the  Its-Lisovyy-Prokhorov conjecture.  Similarly to (\ref{At}) and (\ref{Mt}) we introduce the two spaces 
  \be
\widetilde {\Acal}_L=\Big\{(p,\{t_j\}_{j=1}^N)\;,\; p\in \Acal_L,\; t_j\in \C,\; t_j\neq t_k\Big\}\;,
\la{At1}
\ee 
\be
\widetilde {\Mcal}_L=\Big\{(p,\{t_j\}_{j=1}^N)\;,\; p\in \Mcal_L,\; t_j\in \C,\; t_j\neq t_k\Big\}\;.
\la{Mt1}
\ee
Denote the  pullback of the form $\o_{\Acal}^L$ with respect to the natural projection of   $\At_L$ to $\Acal_L$ by $\wt_{\Acal}^L$  and the  pullback of the form $\o_{\Mcal}$ with respect to the natural projection of $\Mt_L$ to $\Mcal_L$ by $\wt_{\Mcal}^L$.

  \begin{proposition}
  
  The following identity holds between two-forms on $\At_L$:
\be
\widetilde {\Fcal}^* [\wt_{\Mcal}^L]= \wt_{\Acal}^L-\sum_{k=1}^N \d H_k\wedge \d t_k
\la{weak1}
\ee
where $ H_k$ are the   Hamiltonians (\ref{Hkint}).
\end{proposition}
  {\it Proof.} 
 Denote by $2d$ the dimension of the spaces ${\Acal}_L$ and ${\Mcal}_L$.
 Introduce some local Darboux coordinates $(p_i,q_i)$ on ${\Acal}^L$ for the form
 $\o^L_\Acal$ (\ref{wLA}) and also some Darboux coordinates $(P_i,Q_i)$ on ${\Mcal}^L$ for the form $\o^L_\Mcal$ given by (\ref{wLG}).

  We are going to verify (\ref{weak1}) using coordinates $\{t_j\}_{j=1}^N$ and 
  $\{P_j,Q_j\}_{j=1}^d$. Let us split the operator $\d$ into two parts:
  $$
  \d=\d_t + \d_\Mcal
  $$
  where $\d_\Mcal$ is the differential with respect to $\{P_j,Q_j\}_{j=1}^d$.
  Then relation (\ref{HitAM}) can be written as
    \be
\sum_{j=1}^d \d P_j\wedge \d Q_j=\sum_{j=1}^d \d_\Mcal p_i\wedge \d_\Mcal q_i\;.
\la{KSrel}
\ee
The right-hand side can be further rewritten using the Hamilton equations 
 $\f{\p{p_i}}{\p{t_k}}=-\f{\p H_k}{\p q_i}$; $\f{\p{q_i}}{\p{t_k}}=\f{\p H_k}{\p p_i}$ (where the Hamiltonians $H_k$ are given by (\ref{Hkint})). Using
 $$
\d_\Mcal  p_i =\d p_i+ \sum_{k=1}^N \f{\p H_k}{\p q_i} \d t_k\;\hskip0.7cm
\d_\Mcal q_i =\d q_i- \sum_{k=1}^N \f{\p H_k}{\p p_i} \d t_k
$$
one gets
$$
\sum_{i=1}^d \d_\Mcal p_i\wedge \d_\Mcal q_i=\sum_{i=1}^d \d p_i\wedge \d q_i + \sum_{k=1}^N
\d t_k\wedge \sum_{i=1}^d   \left(\f{\p H_k}{\p q_i}  \d q_i +\f{\p H_k}{\p p_i}  \d p_i \right)
$$
\be
-\sum_{\ell<k=1}^{N} \sum_{i=1}^d \left(\f{\p H_\ell}{\p q_i} \f{\p H_k}{\p p_i}-
\f{\p H_\ell}{\p p_i} \f{\p H_k}{\p q_i}\right)\d t_\ell\wedge \d t_k\;.
\la{HHpp}\ee
To simplify the second sum  in (\ref{HHpp}) we recall that
$$
\d H_k= \sum_{i=1}^d   \left(\f{\p H_k}{\p q_i}  \d q_i +\f{\p H_k}{\p p_i}  \d p_i \right) 
+\sum_{\ell=1}^N \f{\p H_k}{\p t_\ell}\Big|_{p_i,q_i=const} \d t_\ell\;;
$$
 thus the second sum can be written as
$$
 \sum_{k=1}^H \d t_k\wedge \d H_k+\sum_{l,k,\, l<k} \left(\f{\p H_k}{\p t_l}\Big|_{p,q}-\f{\p H_l}{\p t_k}\Big|_{p,q}\right) \d t_l\wedge \d t_k\;.
 $$
 Adding all the terms in (\ref{HHpp})  we obtain 
 $$
\sum_{j=1}^{d} \d P_j\wedge \d Q_j = \sum_{i=1}^d \d p_i\wedge \d q_i + 
\sum_{k=1}^N \d t_k \wedge \d H_k - \sum_{\ell<k} \left(\f{\p H_\ell}{\p t_k}\Big|_{p,q}-\f{\p H_k}{\p t_\ell}\Big|_{p,q} + \big\{H_k, H_\ell\big\}\right) \d t_\ell\wedge \d t_k\;.
 $$ 
The  coefficient of $\d t_\ell \wedge \d t_k$ vanishes because the Hamiltonians satisfy the zero--curvature equations implied by commutativity of the flows with respect to $t_j$ and $t_\ell$; in fact in this particular case  they satisfy a stronger compatibility: $\{H_k ,H_\ell\}=0$ and  $\pa_{t_\ell} H_k = \pa_{t_k}H_\ell$.  
Therefore we arrive at 
(\ref{weak1}).
\QED

  Let us show that (\ref{weak1}) implies 
  
  \begin{proposition}[Weak ILP conjecture]
  The following identity holds on the space $\Mt_L$: 
  \be
  \d\Theta_{ILP}^L = 
  \wt_{\Mcal}^L
  \la{wwt}
  \ee
  where 
\be
 \Theta_{ILP}^L=\sum_{j<k}^N\tr A_j A_k \d\log (t_j-t_k)+\sum_{j=1}^N \tr (L_j G_j^{-1} \d_\Mcal G_j)
\la{omweak}\ee
and matrices $G_j$ diagonalizing $A_j$ are chosen to satisfy the Schlesinger equations (\ref{eqGint});
 $\d_\Mcal$ denotes the differential with respect to monodromy coordinates. The form 
 $  \Theta_{ILP}^L  $ 
  is the "weak" version of the form (\ref{ILP}). The form $\wt_{\Mcal}^L $ is the pullback of Alekseev-Malkin form (\ref{wLG}) from $\Mcal_L$ to $\widetilde {\Mcal}_L$.
  \end{proposition}
  {\it Proof.}
  The symplectic potential for the  form $\tilde{\o}_\Acal^L$ 
  can be written as
  \be
  \tilde{\theta}_\Acal^L=\sum_{j=1}^n\tr [L_j G_j^{-1} (\d_t+\d_\Mcal) G_j]\;.
\la{defTKK}
  \ee
We notice that the potential $\tilde{\theta}_\Acal^L$, in contrast to the form $\tilde{\o}_\Acal^L$
itself, is not well-defined on the space  $\At_L$ due to ambiguity $G_j\to G_j D_j$ 
for diagonal $D_j$ in the definition of $G_j$. Under such transformation  
$\theta_\Acal^L$ changes by an exact form.  Therefore  for the purpose of proving (\ref{wwt}) one can pick any concrete representative for each $G_j$. The most natural choice is to assume that $\{G_j\}$  satisfy the system (\ref{eqGint}). 
Then the ``$t$"-part of potential (\ref{defTKK}) can be  computed using \eqref{eqGint} and the definition of the Hamiltonians \eqref{Hkint} to give
\be
\sum_{j=1}^n\tr (L_j G_j^{-1} \d_t G_j) = 2\sum_{j=1}^N H_j \d t_j\;.
\la{tpart}
\ee
 Therefore, the relation (\ref{weak1}) can be rewritten as 
 \be
\widetilde {\Fcal}^* [\wt_{\Mcal}^L]= \d\left(\sum_{k=1}^N \d H_k\wedge \d t_k +\sum_{j=1}^n\tr (L_j G_j^{-1}\d_\Mcal G_j)\right)
\ee
which coincides with (\ref{wwt}).
\QED
 %%%%
  
\paragraph{Comparison of weak and strong ILP conjectures.}
 In spite of the formal similarity, there is a  significant difference between the statements of the weak and strong ILP
 conjectures. In the strong version the form $\sum \tr (L_j \d G_j G_j^{-1})$ is a well-defined form on the phase space $\Acal$
 as well as on its extension $\widetilde {\Acal}$. 
 
 In the weak version the same form is not defined on the space $\Acal^L$ since to get the equality (\ref{wwt}) one needs to take the residues $A_j$ (which are given by a point of $\Acal^L$ up to a conjugation) and then diagonalize  each 
 $A_j$ into $G_j L_j G_j^{-1}$ in a way which is non-local in times $t_j$: the matrices $G_j$'s themselves must satisfy the Schlesinger system (\ref{eqGint}). This requirement can not be satisfied staying entirely within the space $\Acal^L$ and thus $G_j$'s can not be chosen as functionals of $A_j$'s only; their choice encodes  a highly non-trivial $t_j$-dependence which fixes the freedom in the right multiplication of each $G_j$ by a diagonal matrix which also can be time-dependent.
 
 The strong version of the ILP conjecture (Theorem \ref{dTTH})  is a stronger statement since the form $\theta_\Acal$ is  a 1-form defined on the underlying  phase space.

\section{Log--canonical coordinates and symplectic potential}
\la{FGform}
 Here we  summarize results of \cite{BK} where the form $\omega_\Mcal$ was expressed in $\log$-canonical form  an open subspace of highest dimension of $\Mcal$ 
using the (extended) system of
 Fock-Goncharov  coordinates \cite{FG}. This allows to find the corresponding symplectic potential and use it in the definition of the tau-function.

\subsection{Fock-Goncharov coordinates}
\la{secFG}

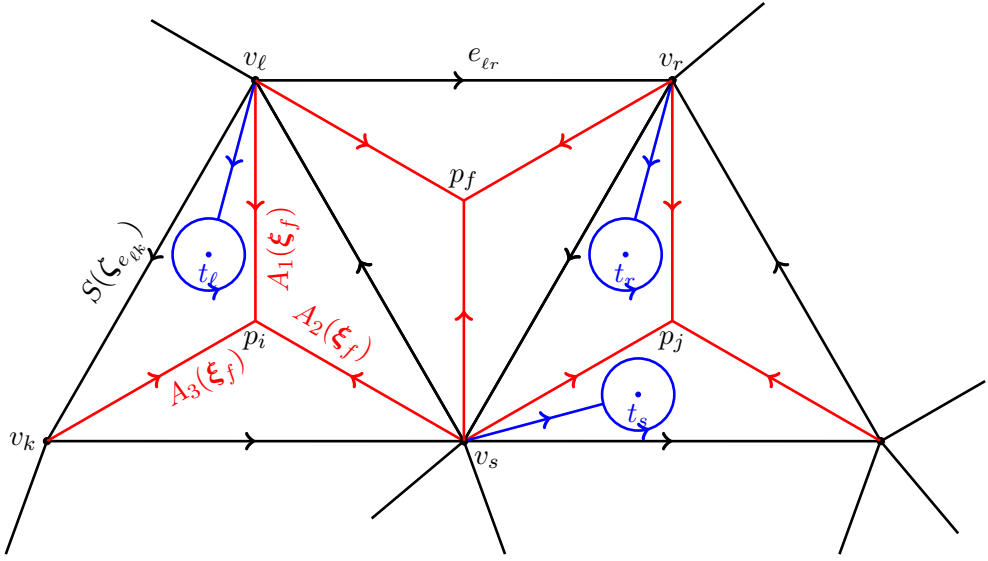
\begin{figure}
    \begin{center}
    \begin{tikzpicture}[scale=1.6]
    \draw [fill] (-30:2) circle[radius=0.03];
    \draw [fill] (90:2) circle[radius=0.03];
    \draw [fill] (210:2) circle[radius=0.03];
    \draw [line width=1pt,  
    postaction={decorate,decoration={markings,mark=at position 0.1666 with {\arrow[line width=1.5pt]{>}}}},
    postaction={decorate,decoration={markings,mark=at position 0.8333 with {\arrow[line width=1.5pt]{>}}}},
    postaction={decorate,decoration={markings,mark=at position 0.5 with {\arrow[line width=1.5pt]{>}}}}](-30:2) to (90:2) to node[sloped, above, pos=0.55] {$S({\bs \zeta}_{e_{_{\ell k}}})$} (210:2) to  cycle; 
    \draw [red,line width=1pt,  postaction={decorate,decoration={markings,mark=at position 0.5 with {\arrow[line width=1.5pt]{<}}}}](0,0) to node[above, pos=0.3, sloped] {$A_2({\bs\xi}_{f})$}(-30:2);
    \draw [red,line width=1pt,  postaction={decorate,decoration={markings,mark=at position 0.5 with {\arrow[line width=1.5pt]{<}}}}](0,0) to node[below, pos=0.3, sloped] {$A_1({\bs\xi}_{f})$}(90:2) coordinate (v1);
    \draw [red,line width=1pt,  postaction={decorate,decoration={markings,mark=at position 0.5 with {\arrow[line width=1.5pt]{<}}}}](0,0) to node[below, pos=0.3, sloped] {$A_3({\bs\xi}_{f})$}(210:2);
    \draw [blue,line width=1pt,  postaction={decorate,decoration={markings,mark=at position 0.5 with {\arrow[line width=1.5pt]{>}}}}]  (v1) to ($ (v1)+ (-105:1.5) $) coordinate(c1);
    \draw [blue, fill=white, line width=1.pt,  postaction={decorate,decoration={markings,mark=at position 0.8 with {\arrow[line width=1.5pt]{>}}}}] (c1) circle [radius=0.3];
    \draw  [blue, fill] (c1) circle[radius=0.02];
    \node [blue] at (c1) [below]{$t_\ell$};
    
    \begin{scope}[shift=(0:3.464106)]
    \draw [fill] (-30:2) circle[radius=0.03];
    \draw [fill] (90:2) circle[radius=0.03];
    \draw [fill] (210:2) circle[radius=0.03];
    \draw [line width=1pt,  
    postaction={decorate,decoration={markings,mark=at position 0.1666 with {\arrow[line width=1.5pt]{>}}}},
    postaction={decorate,decoration={markings,mark=at position 0.8333 with {\arrow[line width=1.5pt]{>}}}}
    ,postaction={decorate,decoration={markings,mark=at position 0.5 with {\arrow[line width=1.5pt]{>}}}}](-30:2) to (90:2) to (210:2) to cycle; 
    \draw [red,line width=1pt,  postaction={decorate,decoration={markings,mark=at position 0.5 with {\arrow[line width=1.5pt]{<}}}}](0,0) to (-30:2);
    \draw [red,line width=1pt,  postaction={decorate,decoration={markings,mark=at position 0.5 with {\arrow[line width=1.5pt]{<}}}}](0,0) to (90:2) coordinate (v1);
    \draw [red,line width=1pt,  postaction={decorate,decoration={markings,mark=at position 0.5 with {\arrow[line width=1.5pt]{<}}}}](0,0) to (210:2) coordinate (v0);
    \draw [blue ,line width=1pt,  postaction={decorate,decoration={markings,mark=at position 0.5 with {\arrow[line width=1.5pt]{>}}}}]  (v1) to ($ (v1)+ (-105:1.5) $) coordinate(c1);
    \draw [blue, fill=white, line width=1.pt,  postaction={decorate,decoration={markings,mark=at position 0.8 with {\arrow[line width=1.5pt]{>}}}}] (c1) circle [radius=0.3];
    \draw  [blue, fill] (c1) circle[radius=0.02];
    \node [blue] at (c1) [below]{$t_r$};
    
    \draw [blue, line width=1pt,  postaction={decorate,decoration={markings,mark=at position 0.5 with {\arrow[line width=1.5pt]{>}}}}]  (v0) to ($ (v0)+ (15:1.5) $) coordinate(c0);
    \node at (v0) [below right] {$v_s$};
    \draw [blue, fill=white, line width=1.pt,  postaction={decorate,decoration={markings,mark=at position 0.8 with {\arrow[line width=1.5pt]{>}}}}] (c0) circle [radius=0.3];
    \draw  [blue, fill] (c0) circle[radius=0.02];
    \node [blue] at (c0) [below]{$t_s$};
    \end{scope}

    \begin{scope}[yscale=-1, shift=(-30:2)]
    \draw [fill] (-30:2) circle[radius=0.03];
    \draw [fill] (90:2) circle[radius=0.03];
    \draw [fill] (210:2) circle[radius=0.03];
    \draw [line width=1pt,  postaction={decorate,decoration={markings,mark=at position 0.1666 with {\arrow[line width=1.5pt]{>}}}},
    postaction={decorate,decoration={markings,mark=at position 0.8333 with {\arrow[line width=1.5pt]{>}}}}
    ,postaction={decorate,decoration={markings,mark=at position 0.5 with {\arrow[line width=1.5pt]{>}}}}](-30:2) to (90:2) to (210:2) to  node[pos=0.55, above]{$e_{_{\ell r}}$} cycle; 
    \draw [red,line width=1pt,  postaction={decorate,decoration={markings,mark=at position 0.5 with {\arrow[line width=1.5pt]{<}}}}](0,0) to (-30:2);
    \draw [red,line width=1pt,  postaction={decorate,decoration={markings,mark=at position 0.5 with {\arrow[line width=1.5pt]{<}}}}](0,0) to (90:2) coordinate (v1);
    \draw [red,line width=1pt,  postaction={decorate,decoration={markings,mark=at position 0.5 with {\arrow[line width=1.5pt]{<}}}}](0,0) to (210:2) coordinate (v0);
    \end{scope}
    
    \draw[ line width = 1pt] (90:2) to ($ (90:2) + (150:1) $);
    \draw[ line width = 1pt] (3.461406,2) to ($ (3.461406,2) + (40:1) $);
    \draw[ line width = 1pt] (-30:2) to ($ (-30:2) + (-70:1) $);
    \draw[ line width = 1pt] (-30:2) to ($ (-30:2) + (-140:1) $);
    \draw[ line width = 1pt] (210:2) to ($ (210:2) + (-110:1) $);
    \draw[ line width = 1pt] ($ (-30:2) + (0:3.461406) $) coordinate(x)  to ($ (x) + (-110:1) $);
    \draw[ line width = 1pt] (x) coordinate(x)  to ($ (x) + (-50:1) $);
    \draw[ line width = 1pt] (x) coordinate(x)  to ($ (x) + (30:1) $);
    \node at (210:2) [left] {$v_k$};
    \node at (90:2) [above] {$v_\ell$};
    \node at (3.461406,2) [above] {$v_r$};
    
    \node at (0,0) [ below] {$p_i$};
    \node at (30:2) [above] {$p_f$};
    \node at ($(30:2) + (-30:2)$) [below] {$p_j$};
    \end{tikzpicture}
    \end{center}
    \caption{The support of the jump matrices $J$. The graph $\Sigma_0$  is in black (the triangulation). }
    \label{figtrian}
    \end{figure}

To define the Fock-Goncharov  coordinates we  introduce the following auxiliary graphs (see Fig. \ref{figtrian}):
\begin{enumerate}
\item
The graph $\Sigma_0$ with $N$ vertices $v_{1},\dots, v_{N}$ which defines a   triangulation of the  $N$-punctured sphere; we assume that each vertex $v_j$ lies in a small neighbourhood of the corresponding pole $t_j$.   Since $\Sigma_0$ is a triangulation there are $2N-4$ faces $\{f_k\}_{k=1}^{2N-4}$  and $3N-6$ edges $\{e_k\}_{k=1}^{3N-6}$; the edges are assumed to be oriented.

\item
Consider a small loop around each $t_k$ (the {\it cherry}) and attach it to  the vertex $v_k$ by an edge (the {\it stem of the cherry}). The cherries are assumed to not intersect the edges of $\Sigma_0$. The union of $\Sigma_0$, the stems and the cherries is denoted by $\Sigma_1$.

The graph $\Sigma_1$ is fixed by $\Sigma_0$ if one chooses  the ciliation 
at each vertex of the graph $\Sigma_0$; the ciliation determines the position of the 
stem of the corresponding cherry.

\item

Choose a point $p_f$ inside each face $f_k$ of $\Sigma_0$  and connect it by edges $\Ec_{f}^{(i)}$, $i=1,2,3$ to the vertices of the face,  oriented towards the point $p_{f}$.  We will denote by $\Sigma$ the graph obtained by the augmentation of $\Sigma_1$ and these new edges. It is the graph $\Sigma$ which will be used to compute the form
$\omega_\Mcal$.

\end{enumerate}

We will make use of the following  notations:
by $\alpha_i$, $i=1,\dots, n-1$ we denote  the simple positive roots of $SL(n)$; by  ${\mathrm h}_i$ the we denote the dual roots: 
\be
\alpha_i:= {\rm diag}( 0,\dots, \!\!\!\mathop{1}^{i-pos}\!\!\!,-1,0,\dots),\qquad
{\mathrm h}_i := \le(
\begin{array}{cc}
(n-i)\1_i & 0\\
0& -i\1_{n-i}
\end{array}
\ri),\qquad 
\tr (\alpha_i {\mathrm h}_k) = n \delta_{ik}\;.
\label{rootprods}
\ee
For  any matrix $M$ we define $M^\star:= \PEM M\PEM$ where $\PEM$ is the ``long permutation" in the Weyl group, 
$$\PEM_{ab}= \delta_{a,n+1-b}\;.$$
In particular 
$\alpha_i^\star = -\alpha_{n-i}$, $ {\mathrm h}_i^\star = -{\mathrm h}_{n-i}\;.$
Let
$$\sigma ={\rm diag} (1,-1,1,-1,\dots)$$
be the signature matrix.

Introduce  the $(n-1)\times(n-1)$ matrix  $\CM$  given by 
\be
\CM_{jk}=\tr ({\mathrm h}_j {\mathrm h}_k)=n^2\left({\rm min}(j,k)-\f{jk}{n}\right)\;.
\la{defCM}
\ee
The matrix $\CM$ coincides with  $n^2 A_{n-1}^{-1}$ with $A_{n-1}$ being the Cartan matrix of $SL(n)$.

The full set of coordinates on $\Mcal$ consists of three groups: the coordinates assigned to vertices of the
graph $\Sigma_0$, to its edges and faces. Below we describe these three groups separately and use them to parametrize the jump matrices of the Riemann-Hilbert problem on the graph $\Sigma$.

\paragraph{ Edge coordinates and jump matrices on $e_j$.}

To each edge $e \in E(\Sigma_0)$ we associate $n-1$ non-vanishing variables 
\be
{\bs z} = \bs z_e =  ( z_1, \dots, z_{n-1} ) \in (\C^\times)^{n-1}
\la{zvar}
\ee
and introduce their exponential counterparts:
\be
{\bf \zeta}={\bf \zeta}_e = (\zeta_{1},\dots, \zeta_{n-1})\in \C^{n-1}\;, \qquad \zeta_j= \f{1}{n}\log z_j^n\;.
\la{vare}\ee 

The jump matrix on the oriented edge $e\in {\bf E}(\Sigma_0)$  is given by 
\be
S({\bf z})  = {\bf z}^{-\boldsymbol h}\PEM \sigma:=  \prod_{j=1}^{n-1} z_j^{-{\mathrm h}_j} \PEM \sigma = { \prod_{\ell=1}^{n-1} z_\ell^{\ell}} \le(
\begin{array}{ccccc}
0 & \dots & & &  (-1)^{n-1} \prod_{j=1}^{n-1}z_{j}^{-n}\\
 & & & \hbox{\reflectbox{$\ddots$}} & 0\\
\vdots &&&\\
0 &  { -z_{n-2}^{-n}}{z_{n-1}^{-n}} &0 & \dots\\
1 & 0 \dots
\end{array}
\ri)
\ee
where $h_i$ are the dual roots  (\ref{rootprods}).
For the inverse matrix we have
$$S^{-1}({\bf z}) =\s \PEM {\bf z}^{\bf h}  = (-1)^{n-1} {\bf z}^{\bf h^\star} \PEM\s \;.$$ 

The notation ${\bf z}^{\bf h}$ stands for 
 \be
{\bf z}^{\bf h} = z_1^{h_1}\dots z_{n-1}^{h_{n-1}}\;.
 \ee 
The sets of variables (\ref{zvar}), (\ref{vare}) corresponding to an oriented edge $e$ of $\Sigma_0$ and the opposite edge $-e$ 
are related as follows:
\be
\label{edgereverse}
{\bs \zeta}_{-e} = (\zeta_{e,n-1} ,\dots, \zeta_{e,1}) 
\;; \qquad 
 {\bf z}_{-e}  :=(-1)^{n-1} (z_{e,n-1}, \dots, z_{e,1})\;.
\ee

\paragraph{Face coordinates and jump matrices on $\Ec_{f}^{(i)}$.}

To each face $f\in F(\Sigma_0)$ (i.e. a triangle of the original triangulation) we associate  
$\frac {(n-1)(n-2)}2$   variables ${\bs \xi}_f = \{\xi_{f;\,abc}:\ \ a, b, c\in \N,\ \ \ a + b + c= n\}$  and their exponential counterparts $x_{f;\,abc}:= {\rm e}^{\xi_{f;\,abc}}$ as follows.

The variables $\xi_{f;\,abc}$ define the jump matrices $A_{i}({\bs\xi}_f)$ on three edges $\{\Ec_{f}^{(i)}\}_{i=1}^3$,  which connect a 
chosen point $p_f$ in each face $f$ of the graph $\Sigma_0$ with its three vertices (these edges are shown in red  in Fig. \ref{figtrian}). The enumeration of vertices $v_1$, $v_2$ and $v_3$  is chosen arbitrarily for each face $f$.
Namely, for  a given vertex $v$ and the face $f$ of $\Sigma_0$ such that $v\in \pa f$ we define the index $f(v) \in \{1,2,3\}$ depending on the enumeration that we have chosen for the three edges $\{\Ec_{f}^{(i)}\}$ lying in the face $f$. For example in Fig \ref{figtrian} for the face $f $ containing point $p_i$  we define  $f(v_\ell)=1$, $f(v_k)=3$ and $f(v_s) = 2$.   

The matrices $A_{1,2,3}(\bs\xi_f)$ are defined following  \cite{FG}. First, the matrix $A_1$ is defined by the formula 
\be
A_1({\bf x}) =\sigma   \le(\prod_{k=n-1}^1 N_k\ri) \PEM \;,
\la{defA1}
\ee 
where  $E_{ik}$ are  the elementary matrices and 
\be
F_i= \1 + E_{i+1,i}\;,\qquad 
H_i(x) := x^{h_i}={\rm diag}(\overbrace{x^{i-n},\dots, x^{i-n}}^{\hbox{$i$ times}}, x^i,\dots x^i)\;  ,\ \ \ \ \ i=1,\dots, n-1\;;
\ee
\be
N_k= \le(\prod_{k\leq  i \leq n-2} 
H_{i+1}(x_{n-i-1, i-k+1, k})
F_i\ri) F_{n-1}\;.
\ee
The matrices $A_2$ and $A_3$ are obtained from $A_1$ by cyclically permuting the indices of the variables: 
\be
A_2({\bs x}) = A_1(\{x_{bca}\})\;,\hskip0.7cm
A_3({\bs x}) = A_1(\{x_{cab}\})\;;
\la{defA23}
\ee
the important  property of  the matrices $A_i$ is the equality 
\be
A_1 A_2 A_3
=\1\;
\la{trivA}
\ee
which guarantees the triviality of total monodromy around the point $p_f$ on each face $f$.
 In the first two non-trivial cases the matrices $A_i$ have the following forms:

\begin{enumerate}
\item[{\bf $SL(2)$:} ]
 there are no  face variables and all matrices $A_i=A$ are given by
\be
\label{SL2A}
A=\left( \ba{cc} 0&1\\ -1&-1
\end{array} \right)\;.
\ee
\item[$SL(3)$: ]
 there is one  parameter $x=x_{111}$ for each face. The  matrices $A_1,A_2$ and $A_3$ coincide in this case, too; they are  given by 
\be
A(x) = \frac 1 x  \left( \begin {array}{ccc} 0&0&1\\ \noalign{\medskip}0&-1&-1
\\ \noalign{\medskip}{x}^{3}&{x}^{3}+1&1\end {array} \right)\; .
\ee
%:

\end{enumerate}

\paragraph{Jump matrices on stems.}

The jump matrix on the stem of the cherry connected to a vertex $v$ is defined from the triviality of  total monodromy around $v$. 

For each vertex $v$ of $\Sigma_0$ of valence $n_v$ the jump matrix on the stem of the  cherry attached to $v$ is given by 
\be
M_v^{0} = \le(\prod_{i=1}^{n_v} A_{f_i}  S_{e_i} \ri)^{-1} 
\la{Mv0}
\ee
where $f_1,\dots, f_{n_v}$ and $e_1,\dots e_{n_v}$ are the faces/edges ordered counterclockwise starting from the stem of the cherry, with the  edges oriented away from the vertex (using if necessary the formula \eqref{edgereverse}).   Since each product $A_{f_i}  S_{e_i}$ is a lower triangular matrix, the matrices $M_v^0$ are also lower--triangular. The diagonal part of $M^0_v$ will be  denoted by $\Lambda_v$ and parametrized as follows:
\be
\Lambda_v
= {\rm diag} \le(m_{v;1} , \frac {m_{v;2}}{m_{v;1}}, \dots, \frac {m_{v;n-1}}{m_{v;n-2}}, \frac 1 { m_{v;n-1}}\ri).
\label{moneig}
\ee  
Notice that the  relation  (\ref{moneig}) can also  be written as  $\Lambda_v={\bf m}_v^{\boldsymbol \a}=\prod_{j=1}^{n-1} m_{v;j}^{\alpha_j}$  where $\alpha_j$ are the roots (\ref{rootprods}).

{
In order to express $\Lambda_v$ in terms of $\zeta$ and $\xi$-coordinates , we enumerate the faces and edges incident at the vertex $v$ by $f_1,\dots, f_{n_v}$ and $e_1,\dots, e_{n_v}$, respectively. We assume the edges to be  oriented away from $v$ using \eqref{edgereverse}. 
We also assume without loss of generality  that the arc $\mathcal E_{f_j}^{(1)}$ is the one connected to the vertex $v$ for all $j=1,\dots, n_v$. Then (see (4.20) of \cite{BK}) we have
\be
\Lambda_v = {\rm e}^{2i\pi L_v}=
 \PEM \left(\prod_{f\perp v} {\bf x}_f^{{\bf h_1}}\right)\left(  \prod_{e\perp v}  {\bf z}_{e}^{{\bf h}} \right) \PEM
\ee
 Introduce now the variables $\mu_{v;\ell}$ via 
\be
\label{muvj}
\mu_{v;n-\ell} = \sum_{f\perp v} \sum_{a+b+c=n\atop a,b,c\geq 1} \xi_{f;abc}\, \CM_{a\ell}
+ \sum_{e\perp v} \sum_{j=1}^{n-1} \zeta_{e;j}\, \CM_{j\ell}
\ee
where the matrix $\CM$ equals to $n^2$ times the inverse Cartan matrix (see (\ref{defCM})). 

The relationship between $\mu$'s and variables $m_v$ is
  $$m_{v;\ell}^n = {\rm e}^{n\mu_{v;\ell}}$$
  i.e. $\mu_{v;\ell}$ defines $m_{v;\ell}$ up to an $n$th root of unity.
Therefore, the entries $\lambda_{v;j}$ of the diagonal matrices $L_v$ are related to  $\mu_{v;j}$ as follows:
\be
\lambda_{v;j}\equiv \f{1}{2\pi i} (\mu_{v;j}-\mu_{v;j-1}) \hskip0.5cm ({\rm mod}\;\; \Z) \;,\hskip0.7cm j=1,\dots,n\;.
\la{lammu}
\ee

}

\paragraph{Vertex coordinates and jump matrices on cherries}

To each vertex $v$ of the graph $\Sigma_0$ we associate a set of $n-1$ non-vanishing complex numbers
$r_{v;i} $, $i=1,\dots,n-1$ in the following way.

\begin{figure}
\begin{center}
\begin{tikzpicture}[scale=02]

\draw[thick,postaction={decorate,decoration={markings,mark=at position 0.9 with {\arrow[black,line width=1.5pt]{>}}}}] circle [radius=0.5];

\draw[postaction={decorate,decoration={{markings,mark=at position 0.15 with {\arrow[black,line width=1.5pt]{<}}}} }]
(0.5,0)--(0: 2);
\draw[fill] (0,0) node[below]{$t$} circle [radius=0.03];
\draw[fill] (0.5,0) node[left]{$\beta$} circle [radius=0.03];
\draw[postaction={decorate,decoration={{markings,mark=at position 0.75 with {\arrow[black,line width=1.5pt]{<}}}} }]
(0.5,0)--(0: 2);
\node at (0.55,0.5){$J_2$};
\node at (0.9, 0.14){$J_1^{-1}$};

\node at (0.45,-0.55){$J_3^{-1}$};
\end{tikzpicture}
\end{center}
\caption{On contribution of one of the loops to the form $d\Theta$}
\label{figstem}
\end{figure}
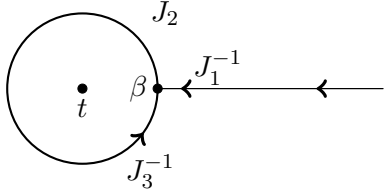

Since the matrix $M^0_v$ is  lower-triangular it can be diagonalized by a lower-triangular matrix $C_v^0$
such that all diagonal entries of $C_v^0$ equal to  1:
\be
M^0_v=C_v^0\Lambda_v (C_v^0)^{-1}\;.
\la{diaM0}
\ee
Any other lower-triangular matrix $C_v$ diagonalizing $M^0_v$ can be written as
\be
C_v= C_v^0 \verB_v
\la{CvCv0}\ee
where the matrix $\verB_v$ (which equals to the diagonal part of $C_v$, $\verB_v=(C_v)^D$),
is parametrized by $n-1$ 
variables $\ver_1, \dots, \ver_{n-1}$ and their logarithmic counterparts  
$$\rho_i = \log r_i \;,\hskip0.7cm i=1,\dots, n-1\;, $$
 as follows (we omit the index $v$ below):
 \be
\verB= \prod_{i=1}^{n-1} \ver_{i}^{ {{\mathrm h}_i}}  = {\bf \ver} ^{\bf h}=
\left(\prod_{i=1}^{n-1} \ver_i^i\right)^{-1} {\rm diag} \Big(\prod_{i=1}^{n-1} \ver_i^n, \,\prod_{i=2}^{n-1} \ver_i^n\;,\, \dots,\ver_{n-2}^n\ver_{n-1}^n,\,  \ver_{n-1}^n,\, 1\Big)\;.
\la{parR}
\ee
 The jump on the boundary of the cherry is defined to be
 \be
 J_v=C_v (z-t_v)^{-L_v}\;.
 \la{jumpcher}
 \ee
{The point of discontinuity of the function $J_v$ on the boundary of the cherry is assumed to coincide with the point where the stem is connected to the cherry (this point is denoted by $\beta$ in Fig.\ref{figstem}).}

\subsection{ Parametrization of  the space $\Mcal$ }

\begin{figure}
\begin{center}
\includegraphics[width=0.4\textwidth]{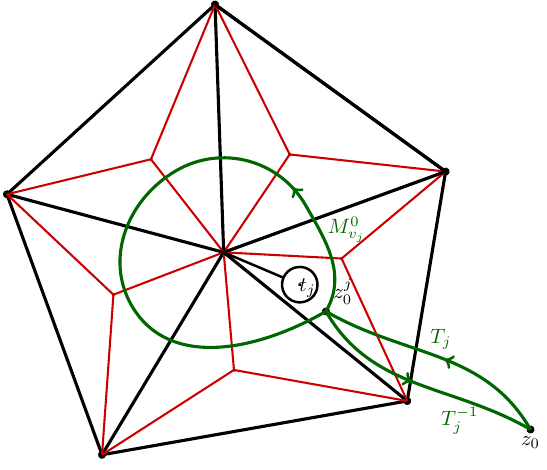}
\end{center}
\caption{ Local triangular monodromy $M_{v_j}^0$ and global monodromy $M_{v_j}^0 = T_j M_{v_j}^0 T_j^{-1}$.
\label{locmonod}
}
\end{figure}

The set of jump matrices on the graph $\Sigma$ constructed in the previous section can be used 
to parametrize the space $\Mcal$. Recall that the vertices of the graph $\Sigma_0$ are in one-to-one correspondence with points $t_j$; thus the vertex connected to the cherry around $t_j$ will be denoted by $v_j$.
To construct the monodromy map as $SL(n)$ representation of $\pi_1( {\bf CP}^1\setminus \{t_1,\dots, t_N\}, \infty)$
{we topologically  identify  the punctured sphere with the complement of connected and simply connected neighbourhoods of the $t_j$'s that contain also the distal vertex of  the stem. The fundamental group of the punctured sphere and of this sphere with deleted neighbourhood is the same. Equivalently, for an element in the fundamental group we choose a representative that does not intersect the cherry and stem. 

Then the  map is then defined as follows; for $\sigma\in  \pi_1( {\bf CP}^1\setminus \{U_1,\dots, U_N\}, \infty)$ the corresponding monodromy is given by 
\be
M_\sigma:=\prod_{e \in \sigma \cap {\bf E}(\Sigma)} J_{e}^{\nu(e,\sigma)}
\ee
where the product is taken in the same order as the order of the edges being crossed by $\sigma$ and $\nu (e,\sigma)\in \{\pm 1\}$ is the orientation of the intersection of the (oriented) edge $e$ and $\sigma$ at the point of intersection. With this definition the analytic continuation of $\Psi$ satisfies $\Psi(z^{\sigma})=\Psi(z) M_\sigma^{-1}$.
}
This allows us to relate  the normalization of the eigenvector matrices $C_j$ with that of the matrices $C_j^0$ \eqref{diaM0}. To this end,  choose $z_0^j$ in the connected region of $\mathbb P^1\setminus \Sigma$ that contains the $j$--th cherry (see Fig. \ref{locmonod}).

Then the monodromy matrix $M_j$ equals to the ordered product of jump matrices at the edges 
of $\Sigma$ crossed by $\sigma_j$ and it
has the form
\be
M_j=T_j M_{v_j}^0 T_j^{-1}
\ee
where the matrix $T_j$ equals to the product of jump matrices on the edges of $\Sigma$ crossed by $\sigma_j$ as it is traversed from $z_0=\infty$ to $z_0^j$.  
Therefore, the diagonal form of the monodromy matrix $M_j$ is:
\be
M_j=C_j \Lambda_j C_j^{-1}, 
\la{MCL}
\qquad 
C_j=T_j C_{v_j}^0 \verB_{v_j} \;.\ee
This determines, unambiguously, the normalization of the matrix $C_j$ in terms of the Fock--Goncharov coordinates, thus providing a complete parametrization of $\wh \Mcal$. 

\subsection{Symplectic form}
The computation of the symplectic form $\d \Theta = \omega_{\M} 
$ is given in  \cite{BK}.

\begin{theorem}[Thm. 4.1 of \cite{BK}]
\label{thmtoric}
In the  coordinate chart parametrized by coordinates
\be
 \bigg\{{\bs z}_e, {\bs x}_f, {\bs r}_v:\ \ e\in E(\Sigma_0), \ f\in F(\Sigma_0), \ v\in V(\Sigma_0)\bigg\}
\ee
the symplectic form $\omega_\M$ \eqref{sympmon}   is given by 
\be\la{ototal}
-2\pi i\;\omega_\M =    \f 1 2  \sum_{v\in V(\Sigma_0)} 
\omega_v+ \f 1 2   \sum_{f\in F(\Sigma_0)} \omega_f +  n \sum_{v\in V(\Sigma_0)} \sum_{i=1}^{n-1}  {  \d \rho_{v;i}    \wedge   \d \mu_{v;i}}\;.
\ee
The variables  $\mu_{v;j}$ are defined by \eqref{muvj}.

The form $\omega_v$ in (\ref{ototal})  is defined as follows:   for each vertex $v\in V(\Sigma_0)$ of valence $n_v$ let $ \{e_1,\dots e_{n_v}\}$ be the incident edges ordered counterclockwise starting from the one on the left of the  stem and oriented away from $v$. Let $\{f_1,\dots, f_{n_v}\}\in F(T)$ be the faces incident to $v$ and counted in counterclockwise order  from the one containing the cherry.  We denote the order relation by $\prec$. Then 
$$
 \omega_v =
 \sum_{e'\prec e \perp v}   \CM_{ij} \d {\zeta}_{e'; i} \wedge    \d  {\zeta}_{e;j}
+
\sum_{f\prec e\perp v} \sum_{a+b+c =n} \sum_{\ell=1}^{n-1} \CM_{f(v),\ell} \d \xi_{f;abc}   \wedge \d \zeta_{e;\ell} 
$$
\be
+
\sum_{e\prec f\perp v} \sum_{a+b+c =n}\sum_{\ell=1}^{n-1}  \CM_{f(v),\ell}  \d \zeta_{e;\ell}  \wedge \d \xi_{f;abc} 
+
\sum_{ f' \prec  f\perp v }  \sum_{a+b+c=n\atop a'\!\!+b'\!\!+c'\!=n}\CM_{f'(v),f(v)}  \d {\xi}_{f'; a'b'c'}  \wedge    \d  {\xi}_{f;abc}
\la{omegav}
\ee
where  the subscript ${f(v)}$ indicates the index $a, b$ or $c$ depending on the value $f(v) \in \{1,2,3\}$, respectively.

The form $\omega_f$ for face $f$ is given by
\bea
\label{omegaf}
 \omega_f
=\sum_{i + j + k = n\atop 
i' + j' + k' = n}
F_{ijk; i'j'k'} \,\, \d \xi_{f; ijk}    \wedge \d \xi_{f; i'j'k'}
\eea
where $F_{ijk; i'j'k'}$ are the following constants
\be
\frac 1 n F_{ijk; i'j'k'}  = 
\Big( k \Delta i -i  \Delta k  \Big) H( \Delta i \Delta k )  + \Big( j \Delta k - k  \Delta j \Big) H( \Delta j \Delta k )  
+
\Big(  i \Delta j -j \Delta i \Big) H( \Delta i \Delta j )  
\la{coefface}\ee
and
$$
 \Delta i = i'-i; \hskip0.7cm \Delta j = j'-j; \hskip0.7cm  \Delta k = k'-k\;;
$$
 $H(x)$ is the Heaviside function:
\be
H(x) = \le\{
\ba{cc}
1  & x>0\\
\frac 1 2  &  x=0\\
0 & x<0
\ea\ri.\;.
\ee
\end{theorem}
\vskip 10pt
We point out that while the coordinates $\bs\xi, \bs\zeta, \bs\rho$ are defined on a covering space of the character variety (with the deck transformations being shifts by integer multiples of $2i\pi$), the symplectic form \eqref{omegav} is defined  on the character variety itself. 
Notice also that for $SL(2)$ and $SL(3)$ the form $\omega_f$ vanishes.

\subsection{Symplectic potential}
\la{SecN}

 We are going to 
 choose a symplectic potential 
 $\theta_\Mcal$ satisfying the equation 
$
\d\theta_\Mcal=\omega_\Mcal
$
for the symplectic form $\omega_\Mcal$
 using the representation (\ref{ototal}). 
For convenience we introduce a uniform notation for coordinates $\zeta_e$ and $\xi_{f;ijk}$; the number of these coordinates equals $\dim \Mcal -(n-1)N$ (we subtract the number of coordinates $\rho_j$ from the total dimension of $\Mcal$). These coordinates we denote collectively by 
$$
\{\kappa_j\}_{j=1}^{\dim \Mcal-(n-1)N}
$$
Then the formula (\ref{ototal}) can be written as
\be
-2\pi i \, \o_\Mcal=\f 1 2 \sum_{j<\ell} n_{j\ell}\d \kappa_j\wedge \d \kappa_\ell+ n\sum_{v\in V(\Sigma_0)} \sum_{j=1}^{n-1}  \d \rho_{v;j}  \wedge  \d \mu_{v;j} 
\la{newco}
\ee
where all $n_{j\ell}$ areinteger numbers and $\mu_{v;j}$ are linear functions of  $\kappa_j$'s.

\begin{definition}\rm
The  symplectic potential  $\theta_\Mcal$ is defined by the following relation:
\be
2 \pi i\, \theta_\Mcal=\f 1 4 \sum_{j<\ell} n_{j\ell}(   \kappa_\ell\d \kappa_j- \kappa_j \d \kappa_\ell  )
+ \f n 2 \sum_{v\in V(\Sigma_0)} \sum_{j=1}^{n-1} (\mu_{v;j} \, \d \rho_{v;j}-  \rho_{v;j}\d  \mu_{v;j})  \;.
\la{symmpoten}
\ee
\end{definition}

Obviously, there exist infinitely many choices of the potential for the form $\omega_\Mcal$.
Our choice (\ref{symmpoten})  is due to Theorem 7.1 of  \cite{BK} which states that the potential $\theta_\Mcal$ (\ref{symmpoten}) remains invariant if any of the cherries is moved to the 
neighbouring face.

\subsection{$SL(2)$ case}
\label{SL2mon}

For $n=2$ the general formula in Thm.  \ref{thmtoric} simplifies considerably to the following (for details see (7.5) of \cite{BK})
\be
-2\pi i \, \omega_\Mcal=\sum_{k=1}^N \le(\sum_{e,e'
\perp v_k\atop e'\prec e}   \d \zeta_{e'}\wedge \d \zeta_{e}
 + 2 \sum_{e\perp v_k} \d\rho_k \wedge \d  \zeta_e\ri)\;;
 \la{omegaMlog}
\ee
the symplectic potential (\ref{symmpoten}) $\theta_\Mcal$  is given by:
\be
2\pi i \theta_\Mcal =\sum_{k=1}^N\le(\frac{1}{2}\sum_{e,e'\perp v_k\atop e' \prec e}( \zeta_{e}  \d \zeta_{e'}-\zeta_{e'}   \, \d \zeta_{e})
 + \sum_{e\perp v_k} (\zeta_e\, \d \rho_k-\rho_k \d \zeta_e)\ri)\;.
 \la{sympotM}
\ee
Notice that the expression (\ref{sympotM}) ``forgets" about the orientation of vertices since the coordinate $\zeta_e$ remains invariant if the orientation of the edge $e$ is changed i.e. when $z_e$ transforms to $-z_e$.
Unlike the form $\omega_\Mcal$, the potential $\theta_\Mcal$ depends on the choice  of  triangulation $\Sigma_0$;     the change of triangulation implies a non-trivial change of $\theta_\Mcal$.

\subsubsection{Change of triangulation. }

One triangulation $\Sigma_0$ can be transformed to any other by a sequence of  ``flips" of the diagonal in
the quadrilateral formed by two triangles with a common edge, see Fig.\ref{flip}.
Let us  assume that the four  cherries attached to the vertices 
are placed as shown in Fig.\ref{flip}.
\begin{figure}
\begin{center}
\includegraphics{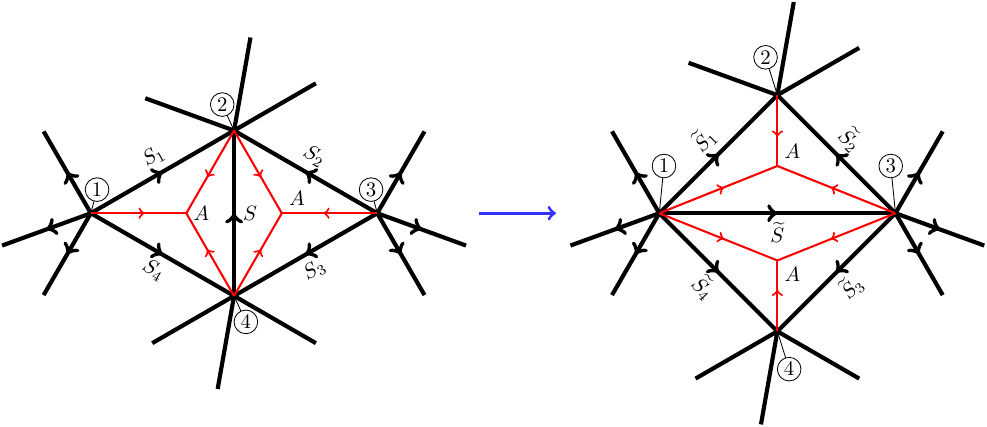}
\end{center}
\caption{
Transformation of edges and jump matrices under an elementary flip of an edge of $\Sigma_0$.}
\la{flip}
\end{figure}
Then, the assumption that all  the monodromies around the four
vertices of these triangles are preserved, implies the following  equations \cite{BK}:
\be
\kat_1=\f{\kappa}{\kappa+1} \kappa_1\;,\hskip0.5cm
\kat_2=(\kappa+1) \kappa_2\;,\hskip0.5cm
\kat_3=\f{\kappa}{\kappa+1} \kappa_3\;,\hskip0.5cm
\kat_4=(\kappa+1) \kappa_4\;,\hskip0.5cm
\kat=\frac{1}{\kappa}
\la{kapkap}
\ee
where $ \kappa_j=z_j^2$, $\kat_j=\tilde{z}_j^2 $, $j=1,\dots,4$; 
  $\tilde\kappa=\tilde{z}^2$ is the variable on the ``flipped" edge.
The variables $r_j$  remain invariant under the change of triangulation due to the choice of cherries positions in Fig.\ref{flip}.

Denote the symplectic potential corresponding to the new triangulation 
 by $\tilde{\theta}_\Mcal$ and introduce the   {\it Rogers dilogarithm}  $L$  which for  $x\geq 0$ is defined by  (see
(1.9) of \cite{Nakanishi}):
\be
L\left(\f{x}{x+1}\right):=  \f{1}{2}\int_0^x \left\{ \f{\log(1+y)}{y}-\f{\log y}{1+y}\right\}\;.
\la{dildef}\ee

As it was shown in Prop.7.1 of \cite{BK},
the symplectic potentials  $\tilde{\theta}_\Mcal$ and $\theta_\Mcal$ are related as follows:
\be
2\pi i (\tilde{\theta}_\Mcal-\theta_\Mcal)= \d \left[L\le(\frac {\ka}{1+\ka} \ri)\right]\;.
\la{ththt}
\ee
Therefore, the function $L\le(\frac {\ka}{1+\ka} \ri)$ is  the generating function of the symplectomorphism corresponding to the elementary flip of the edge of $\Sigma_0$.

\section{ Tau-function as generating function of monodromy symplectomorphism}

\label{tau_fun_mon}

Here we  extend the definition of the  Jimbo-Miwa-Ueno  tau-function by including an explicit description of its dependence on the monodromy data. 

\begin{definition}\la{Deftau}
The tau function is defined by 
 the following set of compatible equations. The equations with respect to $t_j$ are given 
by the Jimbo-Miwa-Ueno formul\ae 
\be
\f{\p \log \tau}{\p t_j}=\f{1}{2}\res{z=t_j}{\rm tr} A^2(z)\;;
\la{taut}
\ee
the equations with respect to coordinates on monodromy manifold $\Mcal$ 
are given by
\be
\d_\Mcal \log \tau=\sum_{j=1}^N \tr (L_j G_j^{-1} \d_\Mcal G_j)-\theta_\Mcal[\Sigma_0]
\la{tauMc}
\ee
where $\theta_\Mcal [\Sigma_0]$ is the symplectic potential (\ref{symmpoten}) for the form $\omega_\Mcal$;
we consider the matrices $G_j$ as (meromorphic) functions on $\widetilde \Mcal$ defined by the formula 
\be
G_j= \Phi(t_j)
\ee
with $\Phi$ the solution of the Riemann-Hilbert problem \eqref{RHPphi}.
\end{definition}

Using  Thm.\ref{ththa} and in particular \eqref{Thetat},   we can rewrite this  definition in an alternative form,
which encodes the complete system (\ref{taut}), (\ref{tauMc}):
\newtheorem{defp}{Definition}
\newenvironment{defprime}[1]
  {\renewcommand{\thedefp}{\ref{#1}$'$}%
   \addtocounter{defp}{-1}%
   \begin{defp}}
  {\end{defp}}

\begin{defprime}{Deftau}
\la{deftaufun}
The  tau-function on $\widetilde {\Mcal}$ is locally defined by equations
\be
\d \log \tau = \Theta -\widetilde \theta_\Mcal[\Sigma_0]
\la{tauMc1}
\ee
where $\Theta$ is the Malgrange form    (\ref{omegaM})   corresponding to solution $\Phi$ \eqref{RHPphi} and $\widetilde \theta_\Mcal[\Sigma_0]$ is the pullback to $\widetilde \Mcal$ of  $\theta_\Mcal[\Sigma_0]$.
\end{defprime}

The formula \eqref{tauMc} means that $\log\tau$  is nothing but  the generating function of the monodromy symplectomorphism: $\d\log\tau$ equals to the difference of symplectic potentials defined in terms of      the (extended) Kirillov--Kostant  symplectic potential $\theta_\Acal$ and the symplectic potential on the monodromy manifold.
It was proven in \cite{Malg} that the  residue of  $\Theta$ along the points of the Malgrange divisor is a positive integer; thus  $\tau$ is actually locally analytic
on $\widetilde{\Mcal}$;  multiplicity of zero of $\tau$ equals to the residue of $\Theta$.

We now analyze in more detail the dependence of $\tau$ on the Fock--Goncharov coordinates. 
The tau-function $\tau$ defined by (\ref{tauMc}) depends on the full set of variables $({\bs z},\,{\bs x},\, {\bs r})$ on
$\Mcal$.
The  right-hand sides  of equations (\ref{tauMc})  depend on the choice of the triangulation $\Sigma_0$  defining the symplectic potential $\theta_\Mcal$.
However, according to Thm.6.1 of \cite{BK}, the potential (\ref{symmpoten}) is independent of the choice of ciliation of the graph
$\Sigma_0$.

The next proposition shows how the tau-function defined in Def.\ref{Deftau} depends on variables $\rho_{j,i}$:
namely, define the second tau-function $\tau_1$ by
\be
\tau_1=\tau \; \exp\left\{- \f{n}{2\pi i} \sum_{j=1}^N\sum_{i=1}^{n-1} \rho_{j;i} \;\mu_{j;i}\right\} \;.
\la{tauNew}
\ee

\begin{proposition}
The tau-function $\tau_1$ (\ref{tauNew}) is independent of the variables $\{r_{j;i}\}, \ j=1\dots N$, $ i=1 \dots n-1$ i.e.
\be
\f{\p \log \tau_1}{\p r_{j;i}}=0\;.
\la{indeptau}
\ee

\end{proposition}
{\it Proof.}
Denote by $G_j^0$ the set of matrices $G_j$ which correspond to all variables $\ver_{j;i}=1$.
Then matrices $G_j$ can be expressed in terms of $G_j^0$ and $\ver_{j;i}$ as follows:
\be
G_j=G_j^0 \verB_j
\la{trG}
\ee
where the diagonal matrix $\verB_j$ is given by (\ref{parR}).
Then,
$$
G_j^{-1} \d  G_j -  (G_j^0)^{-1} \d  G^0_j = R_j^{-1} \d  R_j\;.
$$
Therefore, the first sum in (\ref{tauMc}) gets an additive term equal to
\be
\sum_{j=1}^N \tr\, L_j R_j^{-1} \d R_j\;.
\la{term}\ee
On the other hand, 
matrices $C_j$ transform under (\ref{trG}) in the same way:
\be
C_j=C_j^0 \verB_j
\la{trC}
\ee
where the matrices $C_j^0$ are assumed to be  triangular  with all $1$'s on the diagonal. 

To get the variation of $\theta_\Mcal$ under the transformation (\ref{trC}) we observe that the form $\o_\Mcal$ (\ref{sympmon}) transforms under 
(\ref{trC}) as follows:
$$
\o_\Mcal\to\o_\Mcal +\sum_{j=1}^N \tr \d L_j \wedge R_j^{-1}\d R_j\;.
$$
Therefore, according to our definition of $\theta_\Mcal$, the last sum in this expression should be integrated to give
\be
\theta_\Mcal\to\theta_\Mcal +\sum_{j=1}^N \tr  L_j R_j^{-1}\d R_j
\la{changeth}
\ee
which  cancels against (\ref{term}) (alternatively, one can derive (\ref{changeth}) using the definition of $\mu_{v;i}$ and $\rho_{v;i}$ and (\ref{symmpoten})).

\QED

The  equations for the tau-function  with respect to variables ${\bs z}$ and ${\bs x}$ (or, equivalently, ${\bs \zeta}$ and ${\bs \xi}$)  implied by  Def.\ref{Deftau} can be obtained from  expression 
(\ref{symmpoten}) for the potential $\theta_\Mcal$.  Below we write these equations explicitly in the $SL(2)$ case.

\subsection{$ SL(2)$ tau-function}
\la{SL2tau}
In the $SL(2)$ case the coordinates on 
$\Mcal_N^{SL(2)}$ are given by edge coordinates $\{\zeta_e\}$ and vertex coordinates $\{\rho_k\}_{k=1}^N$;  the potential $\theta_\Mcal$ is given by (\ref{sympotM}). 
Then
$$
L_j=\left(\ba{cc} \lambda_j  & 0 \\ 0 & -\lambda_j\ea\right)=\f{1}{2\pi i} \left(\ba{cc} \mu_j  & 0 \\ 0 & -\mu_j\ea\right)\;.
$$
and the relationship \eqref{tauNew} becomes:
\be
\tau(\{\zeta_e, \rho_j\}, \{t_j\}) =
\tau_1(\{\zeta_e\}, \{ t_j\}) \exp\left\{ \f{1}{\pi i} \sum_{j=1}^N \rho_{j} \;\mu_{j}\right\} .
\ee
where $\mu_j$ is the sum of the $\zeta_e$ for all edges incident to the $j$-th vertex.
  The 
equations  for $\tau_1$ with respect to the edge coordinates take the following form:
\begin{definition}
For a given triangulation $\Sigma_0$ the  tau-function $\tau_1$ of an $SL(2)$ Fuchsian system is defined by the 
system (\ref{taut}) with respect to poles $\{t_j\}_{j=1}^N$ and the following equations with respect to coordinates 
$\{\zeta_{e_j}\}_{j=1}^{3N-6}$:
\be
\f{\p}{\p \zeta_e}\log \tau_1=\sum_{j=1}^N \tr \left(L_j G_j^{-1} \f{\p G_j}{\p \zeta_e}\right)-
\f{1}{4\pi i}\le(\sum_{e'\perp v_1\atop e \prec e'}   \zeta_{e'}-
\sum_{e'\perp v_1\atop e'\prec e}   \zeta_{e'}+\sum_{e'\perp v_2\atop e \prec e'}   \zeta_{e'}-
\sum_{e'\perp v_2\atop e'\prec e}   \zeta_{e'}\ri)
\la{tauze3}
\ee
 where $v_1$ and $v_2$ are vertices of $\Sigma_0$ connected by the edge $e$.
\end{definition}
This definition depends on the choice of triangulation $\Sigma_0$. 
The change of the tau-function $\tau$  under an elementary flip of an edge of the triangulation $\Sigma_0$ acting on the underlying triangulation
follows from (\ref{ththt}):

\begin{proposition}\la{propcoc}
Let $\tau$ and $\tilde{\tau}$ be tau-functions corresponding to triangulations related by the flip 
of the edge $e$
shown in Fig.\ref{flip} .
Then 
\be
\frac{\tilde{\tau}}{\tau}= \exp \left[-\f{1}{2\pi i} L\left(\f{{\rm e}^{2\zeta_e}}{{\rm e}^{2\zeta_e}+1}\right)\right]
\la{cocycle}
\ee
under an appropriate choice of branch of  the Rogers' dilogarithm $L$ (\ref{dildef}).
\end{proposition}

\subsection{Equations with respect to Fock-Goncharov coordinates}

Here we derive equations for $\Psi$, $G_j$ and $\tau$ with respect to an edge coordinate $\zeta$. 

First we notice that for any Riemann-Hilbert problem on an oriented contour $C$ with jump matrix $J$ the variation of the solution of the Riemann-Hilbert problem takes the form :
\be
\delta \Psi \Psi^{-1}(w) = \f{1}{2\pi i} \int_C\f{\Psi_- \delta J J^{-1} \Psi_-^{-1}}{z-w}\d w\;.
\la{genvar}
\ee
The formula (\ref{genvar}) can be easily derived by applying the variation $\delta$ to the equation $\Psi_+=\Psi_- J$ on $C$ which gives $\delta\Psi_+=\delta\Psi_- J+\Psi_- \delta J$
and then solving the resulting non-homogeneous Riemann-Hilbert problem via Cauchy kernel.

We  apply (\ref{genvar}) to variations of  the solution $\Psi$ of  the Riemann-Hilbert problem on the Fock-Goncharov graph with respect to the variable $\zeta$ corresponding to the edge $e$.
Without loss of generality we  assume that the positions of cherries are chosen as in  Fig. \ref{Varform}. 
\begin{figure}[htp]
\begin{center}

 \begin{tikzpicture}[scale=2.5]
 \coordinate (v1) at (-1,1);
 \coordinate (v2) at (1,1);
 \coordinate(v3) at (0,-0.4);

\coordinate (t1) at ($(v1) + (-13:0.7)$);
\coordinate (tt1) at ($(v1) + (-43:0.7)$);
\coordinate (t2) at ($(v2) + (193:0.7)$);
\draw[line width = 1pt, postaction={decorate,decoration={{markings,mark=at position 0.5 with {\arrow[black,line width=1.5pt]{>}}}} }]
(v1) to node[pos=0.45, above] {$e$} (v2);
\draw[line width = 1pt, postaction={decorate,decoration={{markings,mark=at position 0.5 with {\arrow[black,line width=1.5pt]{>}}}} }]
(v2) to node[pos=0.45, above left] {$e_4$} (v3);
\draw[line width = 1pt, postaction={decorate,decoration={{markings,mark=at position 0.5 with {\arrow[black,line width=1.5pt]{>}}}} }]
(v1) to node[pos=0.45, below left] {$e_2$} (v3);

\draw [red!60!black,  line width = 0.6pt, postaction={decorate,decoration={{markings,mark=at position 0.5 with {\arrow[line width=1.5pt]{>}}}} }]   (v1) to (0,0.5);
\draw [red!60!black,  line width = 0.6pt, postaction={decorate,decoration={{markings,mark=at position 0.5 with {\arrow[line width=1.5pt]{>}}}} }]   (v2) to (0,0.5);
\draw [red!60!black,  line width = 0.6pt, postaction={decorate,decoration={{markings,mark=at position 0.5 with {\arrow[line width=1.5pt]{>}}}} }]   (v3) to (0,0.5);

\draw [blue!60!black,  line width = 1pt, postaction={decorate,decoration={{markings,mark=at position 0.5 with {\arrow[line width=1.5pt]{>}}}} }]   (v1) to (t1);
\draw[line width =  1pt, blue!60!black,fill= white, postaction={decorate,decoration={{markings,mark=at position 0.9 with {\arrow[line width=1.5pt]{>}}}} }] (t1) circle[radius=0.13];
\node at (t1) {$t_1$};

\draw [blue!60!black,  line width = 1pt, postaction={decorate,decoration={{markings,mark=at position 0.5 with {\arrow[line width=1.5pt]{>}}}} }]   (v2) to (t2);
\draw[line width =  1pt, blue!60!black,fill= white, postaction={decorate,decoration={{markings,mark=at position 0.9 with {\arrow[line width=1.5pt]{>}}}} }] (t2) circle[radius=0.13];

\draw[line width = 1pt, postaction={decorate,decoration={{markings,mark=at position 0.5 with {\arrow[black,line width=1.5pt]{>}}}} }]
(v1) to node[pos=0.45, above] {$e_1$} ($(v1) + (150:1)$);

\draw[line width = 0.6pt, red!60!black, postaction={decorate,decoration={{markings,mark=at position 0.5 with {\arrow[line width=1.5pt]{>}}}} }]
(v1) to ($(v1) + (70:1)$);
\draw[line width = 0.6pt, red!60!black, postaction={decorate,decoration={{markings,mark=at position 0.5 with {\arrow[line width=1.5pt]{>}}}} }]
(v1) to ($(v1) + (230:1)$);

\draw[line width = 0.6pt, red!60!black, postaction={decorate,decoration={{markings,mark=at position 0.5 with {\arrow[line width=1.5pt]{>}}}} }]
(v2) to ($(v2) + (120:1)$);
\draw[line width = 0.6pt, red!60!black, postaction={decorate,decoration={{markings,mark=at position 0.5 with {\arrow[line width=1.5pt]{>}}}} }]
(v2) to ($(v2) + (-40:1)$);
\draw[line width = 1pt, postaction={decorate,decoration={{markings,mark=at position 0.5 with {\arrow[black,line width=1.5pt]{>}}}} }]
(v2) to node[pos=0.45, above] {$e_3$} ($(v2) + (30:1)$);

\node at (t2) {$t_2$};
\end{tikzpicture}

\end{center}
\caption{ Positioning of the cherries for the computation of the  derivative in $\zeta_e$. } 
\la{Varform}
\end{figure}
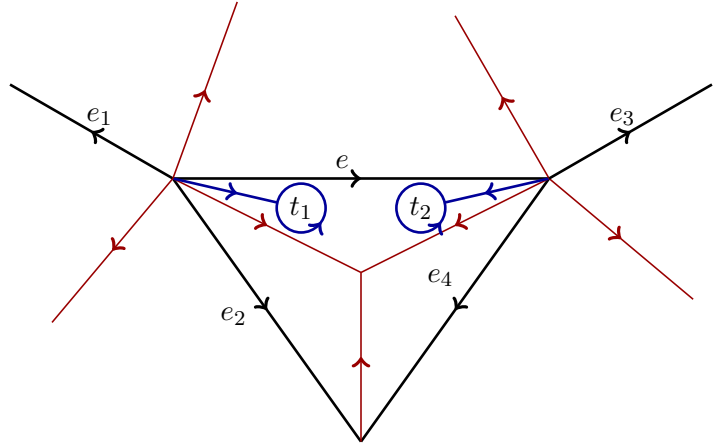

For simplicity we assume that both vertices $v_1$ and $v_2$ connected by $e$ are three-valent but it is not a significant restriction. 

The jump matrices depending on $\zeta$ are the following: the jump matrix on the edge $e$ (and the one on the reverse edge)  according to the general rules in \eqref{edgereverse}  is
\be
S (\zeta)=\left(\ba{cc} 0 & - e^{-\zeta} \\
                       e^{\zeta} & 0 \ea\right)
;\qquad \zeta_{-e} = \zeta+i\pi\;. 
\ee
Furthermore, using the expression \eqref{SL2A} for the matrix $A$  and denoting $S_j= S(\zeta_j)$
the jump matrices $Q_1$ and $Q_2$ on the stems $s_1$ and $s_2$  depend on $\zeta$ as follows: 
$$
Q_1=(S A S_1 A S_2 A)^{-1}=
\left(\ba{cc} e^{\zeta+\zeta_1+\zeta_2} & -e^{-(\zeta+\zeta_1+\zeta_2)} -e^{-\zeta-\zeta_1+\zeta_2} -e^{-\zeta+\zeta_1+\zeta_2}\\
 0   & e^{-(\zeta+\zeta_1+\zeta_2)}\ea\right),
$$
$$
Q_2=(AS_4 AS_3 A S^{-1})^{-1}
= \left(\ba{cc} e^{-(\zeta+\zeta_3+\zeta_4+i\pi )} & 0 \\
e^{\zeta+\zeta_3+\zeta_4+i\pi } +e^{\zeta+\zeta_3-\zeta_4+i\pi } +e^{\zeta-\zeta_3-\zeta_4+i\pi } &
 e^{\zeta+\zeta_3+\zeta_4+i\pi }       \ea\right).
$$
Notice that logarithmic derivatives of the matrices $S$ and $Q_1^{-1}$ and  $Q_2$ with respect  to $\zeta$ are the same and are given by
\be
J_\zeta J^{-1}= -\sigma_3=\left(\ba{cc} -1 & 0\\
                       0 & 1\ea\right)\;.
                       \ee
Then the exponents of monodromy are 
$$\lambda_1=\frac 1{2i\pi}(\zeta+\zeta_1+\zeta_2) \;,\hskip0.7cm \lambda_2=\frac {-1}{2i\pi} (\zeta+\zeta_3+\zeta_4+i\pi)\;.$$
 Thus
\be
L_1=\f{1}{2\pi i} \left(\ba{cc} \zeta+\zeta_1+\zeta_2 &  0 \\
                       0 & -( \zeta+\zeta_1+\zeta_2) \ea\right)\;,\hskip0.7cm
C_1= \left(\ba{cc} 1 &  f_1(\zeta)    \\
                       0 & 1 \ea\right)                       
\la{L1}
\ee
where 
$$
f_1(\zeta)= \f{  e^{-\zeta}(e^{-\zeta_1-\zeta_2} +e^{-\zeta_1+\zeta_2} +e^{\zeta_1+\zeta_2})   }{e^{\zeta+\zeta_1+\zeta_2}- e^{-(\zeta+\zeta_1+\zeta_2)}    }
$$
 and 
\be
L_2=\f{1}{2\pi i} \left(\ba{cc} -(\zeta+\zeta_3+\zeta_4+i\pi ) &  0 \\
                       0 &  \zeta+\zeta_3+\zeta_4+i\pi  \ea\right)\;,\hskip0.7cm
C_2= \left(\ba{cc} 1 &  0    \\
          f_2(\zeta)              & 1 \ea\right)                       
\la{L2}
\ee
where
$$
f_2(\zeta)= \f{  e^{\zeta}(e^{\zeta_3+\zeta_4} +e^{\zeta_3-\zeta_4} +e^{-\zeta_3-\zeta_4} )  }{e^{-(\zeta+\zeta_3+\zeta_4)}- e^{\zeta+\zeta_3+\zeta_4}    }\;.
$$

Introduce the graph $\Sigma_0'$ by identifying the vertex $v_j$ with the corresponding pole $t_j$. The variational formula takes the simplest form if we make an explicit assumption on the growth  \eqref{asint} of the solution $\Psi$ near the poles, that is on the real part of the eigenvalues of $L_j$; indeed it is known \cite{JMU2} that for a given monodromy representation there is a lattice of solutions to the inverse monodromy representation. The reason is simply that the eigenvalues matrices $L_j$ are defined up to addition of  diagonal matrices in $sl_n(\Z)$. For $n=2$ there is therefore a  $\Z^N$ ($N$ being the number of poles) lattice of inverse solutions. The transformations between different solutions in this lattice are referred to as ``discrete Schlesinger transformations''. We  use this observation to shift the real part of $\l_j$'s to a value within the interval $\Re \l_j \in [-\frac 1 2, \frac 12)$.  Excluding the non-generic cases $\Re \l_j = -\frac 1 2$  we have
\begin{theorem}\la{Psizeta} Denote by $\zeta_{jk}$ the FG coordinate corresponding to the 
edge $e_{jk}$. 
Suppose all eigenvalues of matrices $L_j$ satisfy the conditions
\be
|\Re \l_j|<\frac 1 2 \; .
\la{lambda_cond}
\ee
Then the function $\Psi$ satisfies the following system of equations with respect to coordinates $\zeta_{jk}$:
\be
\f{\d \Psi(z)}{\d \zeta_{jk}}= \f{1}{2\pi i} \left[\int_{t_j}^{t_k}\f{\Psi_-(w)\sigma_3\Psi_-^{-1}(w)}{z-w}\d w\right]\Psi(z)\;,\hskip0.7cm z\neq t_j, t_k
\la{varPsi}
\ee
where the integral along the oriented edge $e_{jk}$ of $\Sigma_0'$  in the right hand side is convergent at the endpoints due to condition (\ref{lambda_cond}).
The integral in (\ref{varPsi}) is discontinuous across the  edge $e_{jk}$.
\end{theorem}
{\it Proof.}
We denote $\zeta= \zeta_{jk}$ for simplicity, and $j=1, k=2$. The expression $\pa_{\zeta} J J^{-1}$ is nonzero only on the edge, the two stems and the boundaries of the two cherries. 
A direct computation shows (with the edges oriented as shown in Fig.\ref{Varform}), that the expression of $\pa_{\zeta} J J^{-1}$ on the edge and on the two stems is equal to $\s_3$. 
The jump matrix on the cherry $1$ equals to $J_1=C_1 (z-t_1)^{-L_1}$ and on $c_2$ it equals to $J_2=C_2 (z-t_2)^{-L_2}$.
Then  a direct computation  gives (since ${L_1}_\zeta=\frac {\sigma_3}{2i\pi}$ and ${L_2}_\zeta=\frac{-\sigma_3}{2i\pi} $):
$$
{J_1}_{\zeta} J_1^{-1}={C_1}_{\zeta} C_1^{-1}-\frac{\log(z-t_1)}{2i\pi}C_1 {L_1}_{\zeta} C_1^{-1}=   \left(\ba{cc}0 &  -\f{2 f_1}{1-e^{-4\pi i \lambda_1}} \\
                       0 &0 \ea\right)   - \frac{\log(z-t_1)}{2i\pi} \left(\ba{cc}1 &  -2f_1 \\
                       0 &-1 \ea\right),
$$
$$
{J_2}_{\zeta} J_2^{-1}={C_2}_{\zeta} C_2^{-1}-\frac{\log(z-t_2)}{2i\pi}C_2 {L_2}_{\zeta} C_2^{-1}=   \left(\ba{cc}0 & 0 \\
                        \f{2 f_2}{1-e^{-4\pi i \lambda_2}} &0 \ea\right)   
                       + \frac{\log(z-t_2)}{2i\pi} \left(\ba{cc} -1 &  0 \\
                      -2f_2(\zeta) & 1 \ea\right).
$$
Consider the first cherry  (the second cherry can be treated in parallel); we shall call $\beta$ the point of  intersection of the stem and  the cherry.
Within a neighbourhood containing the cherry, we have $\Psi(z) = \Phi_1(z) (z-t_1)^{L_1} C_1^{-1}$, with $\Phi_1(t_1) = G_1$. 
In the integral \eqref{genvar} the contribution coming from the first cherry is then the integral 
$$
\int_{\beta_+} ^{\beta_-} \frac{\Psi_- \pa_\zeta J_1 J_1^{-1} \Psi_-^{-1}}{w-z} \frac{\d w}{2i\pi} =
$$
\be
\int_{\beta_+} ^{\beta_-} \!\!\!\!\!\!{\Phi_1 \le(f_1 \s_+\le(\frac{i{\rm e}^{2i\pi \l_1}}{\sin(2\pi \l_1)} (w-t_1)^{2\l_1} -\frac{\log(w-t_1)}{i\pi}  \ri)  - \frac{\log(w-t_1)}{2i\pi} C_1^{-1}\s_3 C_1\ri)\Phi_1^{-1}} \frac{\d w}{2i\pi (w-z)} 
\label{deflation}
\ee
where the contour of integration is the circle $|w-t_1|= const$ starting at $\beta_+$ and ending at $\beta_-$ and the branch-cut of the power is the segment from $t_1$ to $\beta$. We also have used that $[C_1,\s_+]=0$. Under the condition $-1< 2\Re \l_1$ the contribution of the integration on the cherry tends to zero in the limit of zero radius of the cherry. 
\QED
\begin{remark}\rm
While the general formula \eqref{deflation} is valid without any restriction on the real parts of $\l_j$'s, the integral on the boundary of the cherries provide some sort of ``regularization'' to the integral along the edge. However, if the  conditions \eqref{lambda_cond} in fulfilled, the regularization are not needed and we arrive at the simplified formula \eqref{varPsi}.
\end{remark}

Introducing the  notation
\be
F(z)= -\Psi_-(z)\sigma_3 \Psi_-^{-1}(z)
\la{defF}
\ee
we can formulate the following
\begin{corollary}\la{Gzeta}
The derivatives of  $G_\ell$ on coordinates $\zeta_{jk}$ take the form:
\be
\f{\d G_\ell}{\d \zeta_{jk}}= \f{1}{2\pi i} \left[\int_{t_j}^{t_k}\f{F(w)}{w-t_\ell}\d w\right] G_\ell\;,\hskip0.7cm \ell\neq j,k
\la{Gz1}\;,
\ee
\be
\f{\d G_j}{\d \zeta_{jk}} G_{j}^{-1} =  \lim_{z\to t_j} \left[\int_{t_j}^{t_k}\f{F(w)}{w-z}\frac{\d w}{2i\pi} +
G_j\le( \frac {i f_j (z-t_j)^{2\l_j} \s_+}{\sin(2\pi \l_j) }  - \frac{\log(z-t_j)}{2i\pi} \s_3\ri)G_j^{-1} \ri]
\la{Gz2},
\ee
\be
\f{\d G_k}{\d \zeta_{jk}} G_{k}^{-1} =  \lim_{z\to t_k} \left[\int_{t_j}^{t_k}\f{F(w)}{w-z}\frac{\d w}{2i\pi} +
G_k\le( \frac {-i f_k (z-t_k)^{-2\l_k} \s_-}{\sin(2\pi \l_k) }  + \frac{\log(z-t_k)}{2i\pi} \s_3\ri)G_k^{-1} \ri]
\la{Gz3}
\ee
where   $f_j = (C_j)_{12}$ and $f_k = (C_k)_{21}$  as defined in  (\ref{L1}) and (\ref{L2}).
\end{corollary}
{\bf Proof.}
The first formula follows from the fact that the connection matrices and exponents at the vertices not connected to the edge are constant under the variation. 
For definiteness assume $j=1$, $k=2$. 
To find $G_1$ ($G_2$ can be treated in the same way) we need to evaluate  $\Phi_1(z) =\Psi(z) C_1 (z-t_1)^{L_1}$ at $z=t_1$.  Differentiating the identity $\Phi_1(z) = \Psi(z) C_1 (z-t_1)^{-L_1}$ with respect to $\zeta$ and taking the limit $z\to t_1$ gives the formula, recalling that $G_1 = \Phi_1(t_1)$.
Indeed 
\be
\pa\Phi_1(z) \Phi_1(z)^{-1} = \pa \Psi(z) \Psi(z)^{-1} + \Phi_1 (z)\le( \frac {i f_1}{\sin (2\pi \l_1)} \s_+ (z-t_1)^{2\l_1} - \frac{\log (z-t_1)}{2i\pi} \s_3 \ri)  \Phi_1(z)^{-1}\;.
\ee
If $\Re 2\l_1>-1$ (which is our standing assumption), in the limit as $z\to t_1$ we can substitute $\Phi_1(z)$ by $G_1$ in the above formula, and we obtain the statement. 

On a related note; the same result is obtained by  looking at the singular behaviour of the integral \eqref{varPsi} (using results on Cauchy--type integrals in \cite{Gakhov}, for example) and simply removing the singular part.
\QED

Now we come to the following
\bp\la{tau1zeta}
The equations (\ref{tauze3}) for the tau-function $\tau_1$ with respect to FG coordinates can be equivalently written as follows:
$$
2\pi i \f{\p}{\p \zeta_{jk}}\log \tau_1=\sum_{\ell=1}^N {\rm reg}\int_{t_j}^{t_k}\f{\tr A_{\ell} F(w)}{w-t_\ell}\d w-\f{1}{2}\left(\sum_{e'\perp v_j\atop e \prec e'}   \zeta_{e'}-
\sum_{e'\perp v_j\atop e'\prec e}   \zeta_{e'}+\sum_{e'\perp v_k\atop e \prec e'}   \zeta_{e'}-
\sum_{e'\perp v_k\atop e'\prec e}   \zeta_{e'}\right)
$$
where $F$ is given by (\ref{defF}); the terms of the sum which  need the regularization correspond to $\ell=j$ and $\ell=k$; they are given by
\be
{\rm reg}\int_{t_j}^{t_k}\f{\tr A_{j} F(w)}{w-t_j}\d w=\lim_{z\to t_j} \left[\int_{t_j}^{t_k}\f{\tr A_j F(w)}{w-z}\d w
 - 2\l_j \log(z-t_j)  \ri]
\ee
and 
\be
{\rm reg}\int_{t_j}^{t_k}\f{\tr A_{k} F(w)}{w-t_k}\d w=\lim_{z\to t_k} \left[\int_{t_j}^{t_k}\f{\tr A_k F(w)}{w-z}\d w
 +2\l_k \log(z-t_k)  \ri]\;.
\ee
\ep
The proof follows from equations (\ref{Gz1}), (\ref{Gz2}) and (\ref{Gz3}).

\QED

 {\bf Acknowledgements.}We thank T.Bridgeland, L.Chekhov, L. Feher, V. Fock, S. Fomin, A. Goncharov, J. Harnad, M. Gekhtman, R. Kashaev, A.Nietzke and M. Shapiro  for  illuminating discussions. We thank also A. Its, O. Lisovyy and A. Prokhorov for 
 clarifying comments. We thank P.Boalch for bringing the references \cite{Jeffrey,Boalch2} to our attention.
Finally, we thank the anonymous referees for helpful commentary. 
The work of M.B. was supported in part by the Natural Sciences and Engineering Research Council of Canada (NSERC) grant RGPIN-2016-06660.
The work of D.K. was supported in part by the NSERC grant
RGPIN/3827-2015.
The completion of this  work was supported by the National Science Foundation under Grant No. DMS-1440140 while the authors were in residence at the Mathematical Sciences Research Institute in Berkeley, California, during the Fall 2019 semester
{\it Holomorphic Differentials in Mathematics and Physics}.

\end{document}